\documentclass[reqno,12pt]{article}
\usepackage{amsmath}
\usepackage{amssymb}
\usepackage{amsthm}
\newtheorem{theorem}{Theorem}

\newtheorem{lemma}{Lemma}

\newtheorem{corollary}[theorem]{Corollary}

\evensidemargin\oddsidemargin

\begin{document}
\title {Estimates of the rate of approximation in the Central Limit Theorem for $L_1$-norm of kernel density estimators.  }

\author{A.Yu. Zaitsev}


\maketitle

\section{Introduction}

\def\cal{\mathcal}

\noindent To fix notation, let $X,$ $X_{1},$ $X_{2},\ldots $ be a sequence
of i.i.d. random variables in ${\bf R}$ with density $f$. Further let $%
\{h_{n}\}_{n\geq 1}$ be a sequence of positive constants such that $%
h_{n}\rightarrow 0$ as $n\rightarrow \infty $. The classical kernel
estimator is defined as
\begin{equation}
f_{n}(x)\stackrel{\rm def}{=}\frac{1}{nh_{n}}\sum_{i=1}^{n}K\left( \frac{%
x-X_{i}}{h_{n}}\right) ,\quad \mbox{for}\;x\in {\bf R},  \label{fnx}
\end{equation}
where $K$ is a kernel satisfying
\begin{equation}
K(u)=0\mbox{,}\quad \mbox{for }|u|>1/2;  \label{k1}
\end{equation}
\begin{equation}
||K||_{\infty }=\sup_{u\in {\bf R}}\left| K(u)\right| =\kappa <\infty ;
\label{k2}
\end{equation}
and
\begin{equation}
\int_{{\bf R}}K(u)\,du=1.  \label{k3}
\end{equation}
Let $||\cdot ||$ denote the $L_{1}({\bf R})$-norm. Write $||K^{2}||=\int_{%
{\bf R}}K^{2}(u)\,du$. For any $t\in {\bf R},$ set
\begin{equation}
\rho (t)=\rho (t,K)\stackrel{\rm def}{=}\frac{\int_{{\bf R}}K(u)\,K(u+t)\,du%
}{||K^{2}||}.  \label{rho}
\end{equation}
Clearly, $\rho (t)$ is a continuous function of $t$, $\left| \rho (t)\right|
\leq 1$, $\rho (0)=1$ and $\rho (t)=0$ for $|t|\geq 1$. Let $Z,$ $Z_{1}$ and
$Z_{2}$ be independent standard normal random variables and set
\begin{equation}
\sigma ^{2}=\sigma ^{2}(K)\stackrel{\rm def}{=}||K^{2}||\int_{-1}^{1}%
\mbox{{\rm cov}}\left( \left| \sqrt{1-\rho ^{2}(t)}\,Z_{1}+\rho
(t)\,Z_{2}\right| ,\left| Z_{2}\right| \right) dt.  \label{sig}
\end{equation}

By definition, any Lebesgue density function $f$ is an element of $L_{1}(%
{\bf R})$. This reason was used by Devroye and Gy\"{o}rfi to justify the
assertion that $||f_{n}-f||$ is the natural distance between a density
function $f$ and its estimator $f_{n}$. In their book, Devroye and
Gy\"{o}rfi \cite{dg}, they posed the question about the asymptotic
distribution of $||f_{n}-f||.$

M. Cs\"{o}rg\H{o} and Horv\'{a}th \cite{ch} were the first who proved a
Central Limit Theorem (CLT) for $||f_{n}-f||_{p},$ the $L_{p}$-norm
distance, $p\geq 1.$ Horv\'{a}th \cite{h} introduced a Poissonization
technique into the study of CLTs for $||f_{n}-f||_{p}.$ The M. Cs\"{o}rg\H{o}
and Horv\'{a}th \cite{ch} and Horv\'{a}th \cite{h} results required some
regularity conditions. Beirlant and Mason \cite{bm} introduced a general
method for deriving the asymptotic normality of the $L_{p}$-norm of
empirical functionals. Mason (see Theorem 8.9 in Eggermont and LaRiccia \cite
{er}) has applied their method to the special case of the $L_{1}$-norm of
the kernel density estimator and proved Theorem \ref{t1.1} below. Gin\'{e},
Mason and Zaitsev \cite{gmz} extended the CLT result of Theorem \ref{t1.1}
to processes indexed by kernels~$K$.

Theorem \ref{t1.1} shows that $\Vert f_{n}-{\bf E\,}f_{n}{\bf \Vert }$ is
asymptotically normal under no assumptions at all on the density $f$.
Centering by ${\bf E\,}f_{n}$ is more natural from a probabilistic point of
view. The estimation of $\Vert f-{\bf E\,}f_{n}{\bf \Vert }$ (if needed) is
a purely analytic problem. The main results of this paper (Theorems \ref
{t1.2}, \ref{t1.3} and \ref{t1.4}) provide estimates of the rate of strong
approximation and bounds for probabilities of moderate deviations in the CLT
of Theorem \ref{t1.1}.\medskip

\begin{theorem}
\label{t1.1} {\it For any Lebesgue density $f$ and for any sequence of
positive constants }$\left\{ h_{n}\right\} _{n\geq 1}$ {\it satisfying $%
h_{n}\rightarrow 0$ and $nh_{n}^{2}\rightarrow \infty $, as $n\rightarrow
\infty ,$ we have}{\rm
\begin{equation}
\frac{\Vert f_{n}-{\bf E\,}f_{n}\Vert -{\bf \,E\,}\Vert f_{n}-{\bf \,E\,}%
f_{n}\Vert }{\sqrt{\mbox{\rm Var}(\Vert f_{n}-{\bf \,E\,}f_{n}\Vert )}}%
\rightarrow _{d}Z  \label{nor}
\end{equation}
}{\it and}
\begin{equation}
\lim_{n\rightarrow \infty }n\,\mbox{{\rm Var}}(\Vert f_{n}-{\bf E\,}%
f_{n}\Vert )=\sigma ^{2}.  \label{con}
\end{equation}
\medskip \end{theorem}

The variance $\sigma ^{2}$ has an alternate representation. Using the
formulas for the absolute moments of a bivariate normal random variable of
Nabeya \cite{n}, we can write
\[
\mbox{{\rm cov}}\left( \left| \sqrt{1-\rho ^{2}(t)}\,Z_{1}+\rho
(t)\,Z_{2}\right| ,\left| Z_{2}\right| \right) =\varphi \left( \rho
(t)\right) ,
\]
where
\begin{equation}
\varphi (\rho )\stackrel{\rm def}{=}\frac{2}{\pi }\left( \rho \arcsin \rho +%
\sqrt{1-\rho ^{2}}-1\right) ,\quad \rho \in \left[ -1,1\right] .
\label{phiofrho}
\end{equation}
It is easy to see that $\varphi (\rho )$ is strictly positive for $\rho \neq
0$. Therefore $\sigma ^{2}>0$. Note that by (\ref{k1}), (\ref{k2}) and (\ref
{sig}),
\begin{equation}
\sigma ^{2}\leq 2\,||K^{2}||\leq 2\,\kappa ^{2}.  \label{sikap}
\end{equation}

In what follows the conditions of Theorem \ref{t1.1} are assumed to hold
unless stated otherwise. We shall denote by $A_{j}$ different universal
constants. We write $A$ for different constants when we do not fix their
numerical values. Throughout the paper, $\theta $ symbolizes any quantity
not exceeding one in absolute value. The indicator function of a set $E$
will be denoted by ${\bf 1}_{E}({\bf \,}\cdot {\bf \,})$. We write $\log
^{*}b=\max \left\{ e,\log b\right\} $.

Let $\eta $ be a Poisson $(n)$ random variable, i.e. a Poisson random
variable with mean $n,$ independent of $X,X_{1},X_{2},\ldots $ and set
\begin{equation}
f_{\eta }(x)\stackrel{\rm def}{=}\frac{1}{nh_{n}}\sum_{i=1}^{\eta }K\left(
\frac{x-X_{i}}{h_{n}}\right) ,  \label{fetax}
\end{equation}
where the empty sum is defined to be zero. Notice that
\begin{equation}
{\bf \,E\,}f_{\eta }(x)={\bf \,E\,}f_{n}(x)=h_{n}^{-1}{\bf \,E\,}K\left(
\frac{x-X}{h_{n}}\right) ,  \label{E}
\end{equation}
\begin{equation}
k_{n}(x)\stackrel{\rm def}{=}n\,\mbox{Var}\left( f_{\eta }(x)\right)
=h_{n}^{-2}{\bf \,E\,}K^{2}\left( \frac{x-X}{h_{n}}\right) ,  \label{kk}
\end{equation}
and
\begin{equation}
n\,\mbox{{\rm Var}}\left( f_{n}(x)\right) =h_{n}^{-2}{\bf \,E\,}K^{2}\left(
\frac{x-X}{h_{n}}\right) -\left\{ h_{n}^{-1}{\bf \,E\,}K\left( \frac{x-X}{%
h_{n}}\right) \right\} ^{2}.  \label{v}
\end{equation}
Define
\begin{equation}
T_{\eta }(x)\stackrel{\rm def}{=}\frac{\sqrt{n}\left\{ f_{\eta }(x)-{\bf %
\,E\,}f_{n}(x)\right\} }{\sqrt{k_{n}(x)}}.  \label{TETTA}
\end{equation}
Let $\eta _{1}$ be a Poisson random variable with mean $1$, independent of $%
X,X_{1},X_{2},\dots $, and set
\begin{equation}
Y_{n}(x)\stackrel{\rm def}{=}\left[ \sum_{j\leq \eta _{1}}K\left( \frac{%
x-X_{j}}{h_{n}}\right) -{\bf E\,}K\left( \frac{x-X}{h_{n}}\right) \right]
\left/ \sqrt{{\bf E\,}K^{2}\left( \frac{x-X}{h_{n}}\right) }\right. .
\label{Yn}
\end{equation}
Let $Y_{n}^{(1)}(x),\dots ,Y_{n}^{(n)}(x)$ be i.i.d. $Y_{n}(x).$ Clearly
(see (\ref{fetax})--(\ref{kk}) and (\ref{TETTA})),
\begin{equation}
T_{\eta }(x)=_{d}\frac{\sum_{i=1}^{n}Y_{n}^{(i)}(x)}{\sqrt{n}}.  \label{Teta}
\end{equation}
Set, for any Borel sets $B,E$,
\begin{equation}
J_{n}(B)\stackrel{\rm def}{=}\sqrt{n}\int_{B}\{|f_{\eta }(x)-{\bf E\,}%
f_{n}(x)|-{\bf E\,}|f_{\eta }(x)-{\bf E\,}f_{n}(x)|\}\,dx,  \label{JnBK}
\end{equation}
\begin{equation}
v_{n}(B,E)\stackrel{\rm def}{=}{\bf E\,}\left[ J_{n}(B){\bf \,}%
J_{n}(E)\right] ,  \label{39a}
\end{equation}
\begin{equation}
\sigma _{n}^{2}(B)\stackrel{\rm def}{=}{\bf E\,}J_{n}^{2}(B)=v_{n}(B,B),
\label{sigk}
\end{equation}
\begin{equation}
{\bf P}(B)\stackrel{\rm def}{=}\int_{B}f(x)\,dx={\bf P}\left\{ X\in
B\right\} ,  \label{Pc}
\end{equation}
and
\begin{equation}
R_{n}(B,E)\stackrel{\rm def}{=}\int_{B}\left( \int_{-1}^{1}\left|
g_{n}(x,t,E)-g(x,t,E)\right| \,dt\right) \,dx,  \label{gnint}
\end{equation}
where
\begin{equation}
g(x,t,E)\stackrel{\rm def}{=}{\bf 1}_{E}(x)\,\mbox{{\rm cov}}\left( \left|
\sqrt{1-\rho ^{2}(t)}\,Z_{1}+\rho (t)\,Z_{2}\right| ,\left| Z_{2}\right|
\right) \,f(x),  \label{gxt}
\end{equation}
\begin{equation}
g_{n}(x,t,E)\stackrel{\rm def}{=}{\bf 1}_{E}(x){\bf 1}_{E}(x+th_{n})\,{\Bbb C%
}_{n}\left( x,x+th_{n}\right) \,\sqrt{f(x)\,f(x+th_{n})},  \label{gnxt}
\end{equation}
\begin{equation}
{\Bbb C}_{n}\left( x,y\right) \stackrel{\rm def}{=}\mbox{{\rm cov}}\left(
\left| \sqrt{1-\rho _{n,x,y}^{2}}Z_{1}+\rho _{n,x,y}\,Z_{2}\right|
,|Z_{2}|\right) ,  \label{cxy}
\end{equation}
$Z_{1}$\ and $Z_{2}$\ are independent standard normal random variables and
\begin{equation}
\rho _{n,x,y}\stackrel{\rm def}{=}{\bf E\,}T_{\eta }(x)\,T_{\eta }(y)={\bf %
\,E\,}Y_{n}(x)\,Y_{n}(y)=\frac{{\bf E\,}\left[ K\left( \frac{x-X}{h_{n}}%
\right) \,K\left( \frac{y-X}{h_{n}}\right) \right] }{\sqrt{{\bf E\,}%
K^{2}\left( \frac{x-X}{h_{n}}\right) {\bf \,E\,}K^{2}\left( \frac{y-X}{h_{n}}%
\right) }}.  \label{ron}
\end{equation}
Note that ${\Bbb C}_{n}\left( x,y\right) $ is non-negative and
\begin{equation}
\sup_{x,y\in {\bf R}}{\Bbb C}_{n}\left( x,y\right) \leq 1.  \label{kn0}
\end{equation}
The following Lemma \ref{l1.1} will be proved in Section 2. It is crucial
for the formulation of the main results of the paper, Theorems \ref{t1.2},
\ref{t1.3} and \ref{t1.4} below.\medskip

\begin{lemma}
\label{l1.1} {\it Whenever }$h_{n}\rightarrow 0${\it \ and }$%
nh_{n}^{2}\rightarrow \infty ${\it , as }$n\rightarrow \infty $,{\it \ there
exist sequences of Borel sets}
\begin{equation}
E_{1}\subset E_{2}\subset \cdots \subset E_{n}\subset \cdots  \label{m7}
\end{equation}
{\it and constants }$\left\{ \beta _{n}\right\} _{n=1}^{\infty }${\it \ and }%
$\left\{ D_{n}\right\} _{n=1}^{\infty }$ {\it such that the density }$f(x)$%
\/ {\it is continuous,\ for} \ $x\in E_{n}$, $n=1,2,\ldots ,${\it \ and
relations}
\begin{equation}
\phi _{n}\stackrel{\rm def}{=}\int_{{\bf R}\backslash
E_{n}}f(x)\,dx\rightarrow 0,\quad \mbox{{\it as} }n\rightarrow \infty ,
\label{psin0}
\end{equation}
\begin{equation}
0<\beta _{n}\stackrel{\rm def}{=}\inf_{y\in E_{n}}f(y)\leq f(x)\leq D_{n}%
\stackrel{\rm def}{=}\sup_{y\in E_{n}}f(y)<\infty ,\quad \mbox{{\it for} }%
x\in E_{n},  \label{bdn}
\end{equation}
{\it and}
\begin{equation}
\varepsilon _{n}\stackrel{\rm def}{=}\sup_{H\in {\cal H}_{0}}\;\sup_{x\in
E_{n}}|f*H_{h_{n}}(x)-I(H)\,f(x)|\rightarrow 0,\quad \mbox{{\it as} }%
n\rightarrow \infty ,  \label{CC}
\end{equation}
{\it are valid, where}
\begin{equation}
I(H)\stackrel{\rm def}{=}\int_{{\bf R}}H(x)\,dx,  \label{IH}
\end{equation}
\begin{equation}
f*H_{h}(x)\stackrel{\rm def}{=}h^{-1}\int_{{\bf R}}f(z)\,H\left( \frac{x-z}{h%
}\right) \,dz,  \label{9*}
\end{equation}
\begin{equation}
{\cal H}_{0}\stackrel{\rm def}{=}\left\{ K,K^{2},\left| K\right| ^{3},{\bf 1}%
\{x:|x|\leq 1/2\}\right\} .  \label{H0}
\end{equation}
{\it Moreover,}
\begin{eqnarray}
\frac{D_{n}^{1/2}}{\beta _{n}^{1/2}}\left( \frac{1}{\left( \beta
_{n}\,nh_{n}\right) ^{1/5}}+\frac{\varepsilon _{n}}{\beta _{n}}\right)
&+&R_{n}(E_{n},E_{n})+\frac{\lambda (E_{n})}{\sqrt{nh_{n}^{2}}}+D_{n}\,h_{n}
\nonumber \\
&+&\frac{D_{n}^{3}\,P_{n}}{\beta _{n}^{3}}+{\Bbb N}_{n}\sqrt{h_{n}}%
\rightarrow 0,\quad \mbox{{\it as} }n\rightarrow \infty ,  \label{tendzero}
\end{eqnarray}
{\it where }$R_{n}(E_{n},E_{n})$ {\it is defined in }{\rm (\ref{gnint}),} $%
\lambda (\,\cdot \,)${\it \ means the Lebesgue measure,}
\begin{equation}
{\Bbb N}_{n}\stackrel{\rm def}{=}\int_{E_{n}}f^{3/2}(x)\,dx,  \label{nan}
\end{equation}
{\it and}
\begin{equation}
P_{n}\stackrel{\rm def}{=}\max_{x\in {\bf R}}{\bf P}\left\{
[x,x+2\,h_{n}]\right\} .  \label{PPPP}
\end{equation}
\medskip \end{lemma}

\begin{theorem}
\label{t1.2} {\it There exists an absolute constant }$A${\it \ such that,
whenever }$h_{n}\rightarrow 0${\it \ and }$nh_{n}^{2}\rightarrow \infty $%
{\it , as }$n\rightarrow \infty $, {\it for any sequence of Borel sets} $%
E_{1},E_{2},\ldots ,E_{n},\ldots ${\it \ satisfying }{\rm (\ref{psin0})--(%
\ref{tendzero})}{\it , there exists an }$n_{0}\in {\bf N}$ {\it such that,}
{\it for any fixed }$x>0${\it \ and for sufficiently large fixed }$n\geq
n_{0}${\it , one can construct on a probability space a sequence of i.i.d.
random variables }$X_{1},X_{2},\ldots $ {\it and a standard normal random
variable }$Z$ {\it such that}
\begin{eqnarray}
&&\hspace{1cm}{\bf P}\left\{ \left| \sqrt{n}{\bf \,}\Vert f_{n}-{\bf E\,}%
f_{n}\Vert -\sqrt{n}{\bf \,E\,}\Vert f_{n}-{\bf E\,}f_{n}\Vert -\sigma
\,Z\right| \geq y_{n}+z+x\right\}  \label{015} \\
&\leq &A\,\Big( \exp \left\{ -A^{-1}\,\sigma ^{-1}x/\tau _{n}^{*}\right\}
+\exp \left\{ -A^{-1}\,\kappa ^{-1}\,\Omega _{n}^{-1/2}z\,\log ^{*}\log
^{*}(z/A\,\kappa \,\Omega _{n}^{1/2})\right\}  \nonumber \\
&+&{\bf P}\left\{ \left| \partial _{n}Z\right| \geq z/2\right\} \Big) ,\quad %
\mbox{{\it for any} }z>0,  \nonumber
\end{eqnarray}
{\it where}
\begin{equation}
\tau _{n}^{*}\stackrel{\rm def}{=}A\,\Psi _{n}^{3/2}\,\left( P_{n}+\psi
_{n}\right) ^{1/2}\rightarrow 0,{\bf \quad }\mbox{{\it as} }n\rightarrow
\infty ,  \label{taun*}
\end{equation}
\begin{equation}
y_{n}\stackrel{\rm def}{=}\frac{A\,\lambda (E_{n})\,\left\| K^{3}\right\| }{%
\left\| K^{2}\right\| \sqrt{nh_{n}^{2}}}+\frac{A\,{\Bbb N}_{n}\sqrt{h_{n}}}{%
\sqrt{\left\| K^{2}\right\| }}\rightarrow 0,{\bf \quad }\mbox{{\it
as} }n\rightarrow \infty ,  \label{yyn}
\end{equation}
\begin{eqnarray}
\partial _{n} &\stackrel{\rm def}{=}&\frac{A\,||K^{2}||}{\sigma \,h_{n}}%
\,\left( {\Bbb L}_{n}+\frac{\varepsilon _{n}\,{\Bbb M}_{n}}{\left\|
K^{2}\right\| }\right)  \label{ddn} \\
&+&A\,\kappa \,\Omega _{n}^{1/2}+\frac{A}{\sigma }\left( \frac{\left\|
K^{3}\right\| \,\lambda (E_{n})}{\left\| K^{2}\right\| \sqrt{nh_{n}^{2}}}%
\right) ^{2}\rightarrow 0,{\bf \quad }\mbox{{\it as} }n\rightarrow \infty ,
\nonumber
\end{eqnarray}
\begin{equation}
{\Bbb L}_{n}\stackrel{\rm def}{=}\int_{E_{n}}\int_{E_{n}}{\bf 1}\{|x-y|\leq
h_{n}\}\,\sqrt{f(x)\,f(y)}\,\,{\Bbb K}_{n}(x,y)\,dx\,dy,  \label{LN}
\end{equation}
\begin{equation}
{\Bbb K}_{n}(x,y)\stackrel{\rm def}{=}\min \left\{ 1-\rho _{n,x,y}^{2},\frac{%
\left\| K^{3}\right\| }{\left( 1-\rho _{n,x,y}^{2}\right) ^{3/2}\left\|
K^{2}\right\| ^{3/2}\sqrt{n\,h_{n}\,f(x)}}\right\}  \label{KN}
\end{equation}
\begin{equation}
{\Bbb M}_{n}\stackrel{\rm def}{=}\int_{E_{n}}\int_{E_{n}}{\bf 1}\{|x-y|\leq
h_{n}\}\,f^{1/2}(x)\,f^{-1/2}(y)\,dx\,dy,  \label{mnab}
\end{equation}
\begin{equation}
\Omega _{n}\stackrel{\rm def}{=}\alpha _{n}+2\,P_{n}+2\,\phi _{n}+\frac{%
4\,||K^{2}||\,R_{n}(E_{n},E_{n})}{\sigma ^{2}}+L(n,{\bf R)}\rightarrow 0,%
{\bf \quad }\mbox{{\it as} }n\rightarrow \infty {\bf ,}  \label{dn3}
\end{equation}
\begin{equation}
\alpha _{n}\stackrel{\rm def}{=}\frac{1296}{5}\left( \tau _{n}^{*}\right)
^{2}\log \frac{1}{\tau _{n}^{*}},  \label{aldef}
\end{equation}
\begin{equation}
\Psi _{n}\stackrel{\rm def}{=}\left\| K^{2}\right\| \,D_{n}\,\beta
_{n}^{-1}\,\kappa ^{2}\,\sigma ^{-4},  \label{PSI}
\end{equation}
\begin{equation}
\psi _{n}\stackrel{\rm def}{=}256\,\kappa ^{2}\,\sigma ^{-2}\,\min \left\{
P_{n},D_{n}\,h_{n}\right\} ,  \label{gamma}
\end{equation}
\begin{equation}
L(n,{\bf R})\stackrel{\rm def}{=}\int_{{\bf R}}\left| h_{n}^{-1}\,{\bf P}%
\{X\in [x-h_{n}/2,x+h_{n}/2]\}-f(x)\right| \,dx\rightarrow 0,\quad %
\mbox{{\it as} }n\rightarrow {\infty }.  \label{lmr}
\end{equation}
\end{theorem}

Denote by $F\{{\,\cdot \,}\}$ ${\bf \,}$and $\Phi \{{\,\cdot \,}\}$ ${\bf \,}
$the probability distributions which correspond to the random variables $%
{\bf \,}\sqrt{n}{\bf \,}\left( \Vert f_{n}-{\bf E}\,f_{n}\Vert -{\bf E\,}%
\Vert f_{n}-{\bf E\,}f_{n}\Vert \right) /\sigma $ and $\,Z$, respectively.
The Prokhorov distance is defined by $\pi (F,\Phi )=\inf \left\{ \varepsilon
:\pi (F,\Phi ,\varepsilon )\leq \varepsilon \right\} $, where
\[
\pi (F,\Phi ,\varepsilon )=\sup_{X}\max \left\{ F\{X\}-\Phi \{X^{\varepsilon
}\}\mbox{, }\Phi \{X\}-F\{X^{\varepsilon }\}\right\} ,\quad \varepsilon >0,
\]
and $X^{\varepsilon }$ is the $\varepsilon $-neighborhood of the Borel set $%
X $.\medskip

\begin{corollary}
\label{c1.1} {\it There exists an absolute constant }$A${\it \ such that,
whenever }$h_{n}\rightarrow 0${\it \ and }$nh_{n}^{2}\rightarrow \infty $%
{\it , as }$n\rightarrow \infty $, {\it for any sequence of Borel sets} $%
E_{1},E_{2},\ldots ,E_{n},\ldots ${\it \ satisfying }{\rm (\ref{psin0})--(%
\ref{tendzero})}{\it , there exists an }$n_{0}\in {\bf N}$ {\it such that,
for sufficiently large fixed $n\geq n_{0}$ and for any }$\varepsilon >0$,
\[
\pi (F,\Phi ,2\,\varepsilon +y_{n}/\sigma )\leq A\,\left( \exp \left\{
-A^{-1}\,\kappa ^{-1}\,\Omega _{n}^{-1/2}\sigma \,\varepsilon \,\log
^{*}\log ^{*}(\sigma \,\varepsilon /A\,\kappa \,\Omega _{n}^{1/2})\right\}
\right.
\]
\[
\left. +\;\exp \left\{ -A^{-1}\,\varepsilon /\tau _{n}^{*}\right\} +{\bf P}%
\left\{ \left| \partial _{n}Z\right| \geq \sigma \,\varepsilon /2\right\}
\right)
\]
{\it and}
\[
{\Bbb \pi }(F,\Phi )\leq y_{n}/\sigma +A\,\tau _{n}^{*}\,\log ^{*}\left(
1/\tau _{n}^{*}\right)
\]
\[
+\,A\,\kappa \,\Omega _{n}^{1/2}\sigma ^{-1}\,\log ^{*}\left( \sigma /\kappa
\,\Omega _{n}^{1/2}\right) /\log ^{*}\log ^{*}(\sigma /\kappa \,\Omega
_{n}^{1/2})+A\,\partial _{n}\,\sigma ^{-1}\,\sqrt{\log ^{*}\left( \sigma
/\partial _{n}\right) },
\]
{\it where }$\tau _{n}^{*},y_{n},\Omega _{n},\partial _{n}${\it \
are defined in} {\rm (\ref{taun*})--(\ref{lmr}).}\medskip
\end{corollary}

\begin{theorem}
\label{t1.3} {\it There exists an absolute constant }$A${\it \ such that,
whenever }$h_{n}\rightarrow 0${\it \ and }$nh_{n}^{2}\rightarrow \infty $%
{\it , as }$n\rightarrow \infty $, {\it for any sequence of Borel sets} $%
E_{1},E_{2},\ldots ,E_{n},\ldots ${\it \ satisfying }{\rm (\ref{psin0})--(%
\ref{tendzero})},{\it \ there exists an }$n_{0}\in {\bf N}$ {\it such that,\
for sufficiently large fixed }$n\geq n_{0}${\it \ and} {\it for any fixed }$%
b $\ {\it satisfying }$\tau _{n}^{*}\leq A^{-1}b${\it , }$b\leq 1${\it , one
can construct on a probability space a sequence of i.i.d. random variables }$%
X_{1},X_{2},\ldots $ {\it and a standard normal random variable }$Z${\it \
such that}
\begin{equation}
\hspace{-1cm}{\bf P}\left\{ \left| \sqrt{n}{\bf \,}\Vert f_{n}-{\bf E\,}%
f_{n}\Vert -\sqrt{n}{\bf \,E\,}\Vert f_{n}-{\bf E\,}f_{n}\Vert -\sigma
\,Z\right| \right.  \label{qwer}
\end{equation}
\begin{equation}
\hspace{1cm}\geq \left. A\,\sigma \,\exp \{-b^{2}/72\left( \tau
_{n}^{*}\right) ^{2}\}+y_{n}+z+x\right\}  \nonumber
\end{equation}
\begin{equation}
\leq A\,\Big( \exp \left\{ -A^{-1}\,\sigma ^{-1}x/\tau _{n}^{*}\right\}
+\exp \{-A^{-1}\,\kappa ^{-1}\,\Omega _{n}^{-1/2}z\,\log ^{*}\log
^{*}(z/A\,\kappa \,\Omega _{n}^{1/2})\}  \nonumber
\end{equation}
\begin{equation}
{}+{\bf P}\left\{ b\,\left| Z\right| >A^{-1}\,\sigma ^{-1}x\right\} +{\bf P}%
\left\{ \left| \partial _{n}Z\right| \geq z/2\right\} \Big) ,\quad %
\mbox{{\it for any} }x,z>0,  \nonumber
\end{equation}
{\it where }$\tau _{n}^{*},y_{n},\Omega _{n},\partial _{n}${\it \
are defined in} {\rm (\ref{taun*})--(\ref{lmr})}.\medskip
\end{theorem}

In the formulations of Theorems \ref{t1.2} and \ref{t1.3} and Corollary \ref
{c1.1}, the numbers $n_{0}$ depend on $\{h_{n}\}_{n\geq 1}$, $%
\{E_{n}\}_{n\geq 1}$, $f$ and $K$.
\medskip

Comparing Theorems \ref{t1.2} and \ref{t1.3}, we observe that in Theorem \ref
{t1.2} the probability space depends essentially on $x$, while in the
statement of Theorem \ref{t1.3} inequality (\ref{qwer}) is valid on the same
probability space (depending on $b$) for any $x>0$. However, (\ref{qwer}) is
weaker than (\ref{015}) for some values of $x$. The same rate of
approximation (as in (\ref{015})) is contained in (\ref{qwer}) if $b^{2}\geq
72\left( \tau _{n}^{*}\right) ^{2}\log (1/\tau _{n}^{*})$ and $x\geq
b^{2}\sigma /\tau _{n}^{*}$ only.
Denote now by $F({\bf \,\cdot \,})$\ and $\Phi ({\bf \,\cdot \,})$\ the
distribution functions of the random variables $\sqrt{n}{\bf \,}\left( \Vert
f_{n}-{\bf E\,}f_{n}\Vert -{\bf E\,}\Vert f_{n}-{\bf E\,}f_{n}\Vert \right)
/\sigma $ and $\,Z$,\ respectively. For example, $\Phi (x)=\Phi \left\{
(-\infty ,x]\right\} $. The following statement about moderate deviations
follows from Theorem~\ref{t1.2}.\medskip

\begin{theorem}
\label{t1.4} {\it Under the conditions of Theorem }$\ref{t1.2}${\it , we
have }$$F(-x)/\Phi (-x)\rightarrow 1\quad\hbox{\it \ and }\quad\left( 1-F(x)\right)
/\left( 1-\Phi (x)\right) \rightarrow 1\quad\hbox{\it  as }n\rightarrow \infty, $$
{\it  if }
\[
0<x=x_{n}=o\left( \min \left\{ \left( \tau _{n}^{*}\right) ^{-1/3},{\bf %
\,\,\,}\Omega _{n}^{-1/6}\left( \log ^{*}\log ^{*}(1/\,\Omega _{n})\right)
^{1/3},{\bf \,\,\,}y_{n}^{-1},{\bf \,\,\,}\partial _{n}^{-1/2}\right\}
\right) .
\]
\end{theorem}

The choice of sets $E_{n}$, which are involved in the formulations of our
results, is not unique. Lemma \ref{l1.1} ensures that, for {\it any} density
$f$, there exist sets $E_{n}$ such that the quantities $\tau
_{n}^{*},y_{n},\Omega _{n}$ and $\partial _{n}$ tend to zero. The
optimization of the choice of $E_{n}$ is a separate problem. However, for
sufficiently regular densities $f$, it is not difficult to choose $E_{n}$ so
that the rate of approximation is good enough, see the examples below. In
our treatment of these examples, we shall use the fact that the function $%
\varphi (\rho )$ in (\ref{phiofrho}) satisfies the Lipshitz condition $%
\left| \varphi (\rho _{1})-\varphi (\rho _{2})\right| \leq \left| \rho
_{1}-\rho _{2}\right| $.\medskip

\noindent {\it Example} 1. Consider the density $f$ of the form $%
f(x)=\sum_{j=1}^{m}r_{j}(x){\bf \,1}_{{\cal J}_{j}}(x)$, where functions $%
r_{j}({\bf \,}\cdot {\bf \,})>0$ satisfy the Lipshitz condition
\[
\left| {\bf \,}r_{j}(x)-r_{j}(y)\right| \leq C\left| x-y\right| ^{\gamma
},\quad 0<\gamma \leq 1,\quad \quad \mbox{for }x,y\in {\cal J}_{j},\quad
j=1,2,\ldots ,m,
\]
where constants $C$ and $\gamma $ are independent of $j$ and ${\cal J}%
_{j}=[a_{j},b_{j})$, $a_{j}<b_{j}$, $j=1,2,\ldots ,m$, is a finite
collection of disjoint intervals. Assume that the values of functions $r_{j}$
are separated from zero and infinity:
\[
0<\beta \leq r_{j}(x)\leq D<\infty \quad \quad \mbox{for }x\in {\cal J}%
_{j},\quad j=1,2,\ldots ,m.
\]
Choose
\[
E_{n}=\bigcup_{j=1}^{m}[a_{j}+h_{n}/2,b_{j}-h_{n}/2].
\]
Without loss of generality we assume $a_{j}+h_{n}/2<b_{j}-h_{n}/2$ and $%
h_{n}\leq 1/4$. Then it is easy to estimate $\phi _{n}=O\left( h_{n}\right) $%
, $\beta \leq \beta _{n}\leq D_{n}\leq D$, $\varepsilon _{n}=O\left(
h_{n}^{\gamma }\right) $, $P_{n}=O\left( h_{n}\right) $, $\Psi _{n}=O\left(
1\right) $, $\psi _{n}=O\left( h_{n}\right) $, $\lambda (E_{n})=O\left(
1\right) $, ${\Bbb N}_{n}=O\left( 1\right) $, $y_{n}=O\left( 1/\sqrt{%
nh_{n}^{2}}+\sqrt{h_{n}}\right) $, $L(n,{\bf R})=O\left( h_{n}^{\gamma
}\right) $, $\tau _{n}^{*}=O\left( \sqrt{h_{n}}\right) $, $\alpha
_{n}=O\left( h_{n}\log \frac{1}{h_{n}}\right) $, $R_{n}(E_{n},E_{n})=O\left(
h_{n}^{\gamma }\right) $, $\Omega _{n}=O\left( h_{n}\log \frac{1}{h_{n}}%
+h_{n}^{\gamma }\right) $, ${\Bbb L}_{n}=O\left( h_{n}\left( nh_{n}\right)
^{-1/5}\right) $, ${\Bbb M}_{n}=O\left( h_{n}\right) $,
\[
\partial _{n}=O\left( \sqrt{h_{n}\log \frac{1}{h_{n}}}+h_{n}^{\gamma
/2}+\left( nh_{n}\right) ^{-1/5}+\frac{1}{nh_{n}^{2}}\right) .
\]
Thus, the statement of Theorem \ref{t1.4} is valid for
\begin{eqnarray*}
0 &<&x=x_{n}=o\left( \min \left\{ {\bf \,}h_{n}^{-1/6}\left( \log \frac{1}{%
h_{n}}\right) ^{-1/6}\left( \log \log \frac{1}{h_{n}}\right) ^{1/3},\right.
\right. \\
&&\left. \left. h_{n}^{-\gamma /6}\left( \log \log \frac{1}{h_{n}}\right)
^{1/3},\,\,\,\,\left( nh_{n}\right) ^{1/10},\,\,\,\,\left( nh_{n}^{2}\right)
^{1/2}{\bf \,}\right\} \right) .
\end{eqnarray*}
\medskip

\noindent {\it Example} 2. Consider the standard normal density $%
f(x)=e^{-x^{2}/2}/\sqrt{2\pi }$. Choose
\[
E_{n}=\left[ -\sqrt{2^{-1}\,\log \frac{1}{h_{n}}},\sqrt{2^{-1}\,\log \frac{1%
}{h_{n}}}\right] .
\]
Without loss of generality we assume $h_{n}\leq 1/4$. Then $\phi
_{n}=O\left( h_{n}^{1/4}\right) $, $\beta _{n}^{-1}=O\left(
h_{n}^{-1/4}\right) $, $D_{n}=O\left( 1\right) $, $\varepsilon _{n}=O\left(
h_{n}\right) $, $P_{n}=O\left( h_{n}\right) $, $\Psi _{n}=O\left(
h_{n}^{-1/4}\right) $, $\psi _{n}=O\left( h_{n}\right) $, $L(n,{\bf R}%
)=O\left( h_{n}\right) $, $\tau _{n}^{*}=O\left( h_{n}^{1/8}\right) $, $%
\alpha _{n}=O\left( h_{n}^{1/4}\log \frac{1}{h_{n}}\right) $, $%
R_{n}(E_{n},E_{n})=O\left( h_{n}\right) $, $\Omega _{n}=O\left(
h_{n}^{1/4}\log \frac{1}{h_{n}}\right) $, ${\Bbb L}_{n}=O\left( h_{n}\left(
nh_{n}\right) ^{-1/5}\right) $, ${\Bbb M}_{n}=O\left( h_{n}\sqrt{\log \frac{1%
}{h_{n}}}\right) $, ${\Bbb N}_{n}=O\left( 1\right) $, $\lambda
(E_{n})=O\left( \sqrt{\log \frac{1}{h_{n}}}\right) $, $$y_{n}=O\left( \sqrt{%
\log \frac{1}{h_{n}}}/\sqrt{nh_{n}^{2}}+\sqrt{h_{n}}\right) ,$$ $$\partial
_{n}=O\left( h_{n}^{1/8}\sqrt{\log \frac{1}{h_{n}}}+\left( nh_{n}\right)
^{-1/5}+\frac{\log \frac{1}{h_{n}}}{nh_{n}^{2}}\right) .$$ The statement of
Theorem \ref{t1.4} is valid for
\begin{equation*}
0<x=x_{n}=o\left( \min {\bf \,}\left\{ h_{n}^{-1/24}\left( \log \frac{1}{%
h_{n}}\right) ^{-1/6}\left( \log \log \frac{1}{h_{n}}\right) ^{1/3}, \quad\left( nh_{n}^{2}\right) ^{1/2}\left( \log \frac{%
1}{h_{n}}\right) ^{-1/2}{\bf \,}\right\} \right) .
\end{equation*}
\medskip

\noindent {\it Example} 3. Consider the density
\[
f(x)=f_{\gamma }(x)=\left\{
\begin{tabular}{l}
$\left| x\right| ^{-\gamma }(1-\gamma ),\quad $ $\,$ $0<x\leq 1,$ \\
$\quad \quad 0,\quad \quad \quad \quad \quad $otherwise,
\end{tabular}
\right. \quad 0<\gamma <1.
\]
Choose $\alpha =\frac{1-\gamma }{1+2\gamma }$ and $E_{n}=\left[
h_{n}^{\alpha },1-h_{n}\right] $. Without loss of generality we assume $%
h_{n}\leq 1/8$. Then it is easy to estimate $\phi _{n}=O\left(
h_{n}^{(1-\gamma )\alpha }\right) $, $\beta _{n}^{-1}=O\left( 1\right) $, $%
D_{n}=O\left( h_{n}^{-\gamma \alpha }\right) $, $\varepsilon _{n}=O\left(
h_{n}^{1-(1+\gamma )\alpha }\right) $, $P_{n}=O\left( h_{n}^{1-\gamma
}\right) $, $\Psi _{n}=O\left( h_{n}^{-\gamma \alpha }\right) $, $\psi
_{n}=O\left( h_{n}^{1-\gamma \alpha }\right) $, ${\Bbb N}_{n}=O\left(
h_{n}^{(1-3\gamma /2)\alpha }+{\bf 1}\left\{ \gamma =2/3\right\} \log \frac{1%
}{h_{n}}\right) $, ${\Bbb L}_{n}=O\left( h_{n}\left( nh_{n}\right)
^{-1/5}\right) $, ${\Bbb M}_{n}=O\left( h_{n}\right) $, $\lambda
(E_{n})=O\left( 1\right) $, $y_{n}=O\left( 1/\sqrt{nh_{n}^{2}}%
+h_{n}^{(1-3\gamma /2)\alpha }\sqrt{h_{n}}\right) $, $%
R_{n}(E_{n},E_{n})=O(h_{n}^{1-2\gamma \alpha })$, $\tau
_{n}^{*}=O(h_{n}^{(1-\gamma -3\gamma \alpha )/2})$, $\alpha _{n}=O\left(
h_{n}^{1-\gamma -3\gamma \alpha }\log \frac{1}{h_{n}}\right) $, $L(n,{\bf R}%
)=O\left( h_{n}^{1-\gamma }\right) $, $\partial _{n}=O\left( \Omega
_{n}^{1/2}+\left( nh_{n}\right) ^{-1/5}+\frac{1}{nh_{n}^{2}}\right) $, $%
\Omega _{n}=O\left( h_{n}^{1-\gamma -3\gamma \alpha }\log \frac{1}{h_{n}}%
+h_{n}^{(1-\gamma )\alpha }\right) $. The statement of Theorem \ref{t1.4} is
valid for
\begin{eqnarray}
\hspace{0.7cm}0 &<&x_{n}=o\left( \min \left\{ {\bf \,}h_{n}^{-(1-\gamma
)^{2}/6(1+2\gamma )}\left( \log \frac{1}{h_{n}}\right) ^{-1/6}\!\!\left(
\log \log \frac{1}{h_{n}}\right) ^{1/3}\!\!,\right. \right.   \label{llog} \\
&&\left. \hspace{3cm}\left. \left( nh_{n}\right) ^{1/10},\,\,\,\,\left(
nh_{n}^{2}\right) ^{1/2}{\bf \,}\right\} \right) .  \nonumber
\end{eqnarray}
Note that the logarithmic factor in (\ref{llog}) could be slightly improved
by means of a more careful choice of the intervals $E_{n}$.$\medskip $

When estimating ${\Bbb L}_{n}$ in the examples, we used the fact that, by (%
\ref{bdn}) and (\ref{KN}), for $x,y\in E_{n}$, we have
\begin{equation}
{\Bbb K}_{n}(x,y)\leq \frac{A\,\left\| K^{3}\right\| ^{2/5}}{\left\|
K^{2}\right\| ^{3/5}\left( nh_{n}\,f(x)\right) ^{1/5}}\leq \frac{A\,\left\|
K^{3}\right\| ^{2/5}}{\left\| K^{2}\right\| ^{3/5}\left( \beta
_{n}\,nh_{n}\right) ^{1/5}}.  \label{KN0}
\end{equation}
For some densities $f$ and kernels $K$, formula (\ref{LN}) may give sharper
bounds. For example, if $K=f={\bf 1}\{x:|x|\leq 1/2\}$ and $%
E_{n}=[h_{n}/2,1-h_{n}/2]$, one can show that ${\Bbb L}_{n}=O\left(
h_{n}\left( nh_{n}\right) ^{-2/5}\right) $. This is better than the rates
given in Examples 1 and 3.

Studying the examples and analyzing the statements of Theorems \ref{t1.2},
\ref{t1.3} and \ref{t1.4}, we see that the rates of normal approximation
become worse when the density $f$ is non-smooth or has too small or too
large values. To show that this is essential, let us consider a scheme of
series, where the density $f$ may be depending on $n$. Namely, let
\[
f(x)=\left( 2a_{n}^{-1}\right) {\bf 1}_{[-a_{n},a_{n}]}(x),
\]
where $a_{n}$ may tend to zero or to infinity as $n\rightarrow \infty $. It
is not difficult to understand that we can choose $a_{n}$ tending to
infinity so fast that, with probability tending to 1, the intervals $\left[
X_{i}-h_{n}/2,X_{i}-h_{n}/2\right] $, $i=1,2,\ldots ,n$, are disjoint and
the distribution of $\sqrt{n}{\bf \,}\Vert f_{n}-{\bf E\,}f_{n}\Vert -\sqrt{n%
}{\bf \,}\left( \Vert K\Vert +1\right) $ converges to the degenerate
distribution ${\Bbb E}_{0}$ concentrated at zero. On the other hand, we can
choose $a_{n}$ tending to zero so fast that $\sqrt{n}{\bf \,}\Vert f_{n}-%
{\bf E\,}f_{n}\Vert $ converges to the same degenerate distribution ${\Bbb E}%
_{0}$ since it behaves as in the case where ${\bf P}\left\{ X=0\right\} =1$.
Thus, if all non-zero values of $f$ are very large or very small, then the
distribution of $\sqrt{n}{\bf \,}\left( \Vert f_{n}-{\bf E\,}f_{n}\Vert -%
{\bf E\,}\Vert f_{n}-{\bf E\,}f_{n}\Vert \right) $ is far from that of $%
\sigma {\bf \,}Z$.

Sections 2--5 are devoted to the proof of Theorems \ref{t1.2}, \ref{t1.3}
and \ref{t1.4}. In the proof, we shall use the Poissonization of the sample
size, considering integrals $\int_{E_{n}}\{|f_{\eta }-{\bf E\,}f_{n}|-{\bf %
E\,}|f_{\eta }-{\bf E\,}f_{n}|\}$ instead of $\int_{E_{n}}\{|f_{n}-{\bf E\,}%
f_{n}|-{\bf E\,}|f_{n}-{\bf E\,}f_{n}|\}$. This allows us to use
independence properties of the Poisson point process $\left\{ X_{1},\ldots
,X_{\eta }\right\} $. In Section 2, we prove Lemma \ref{l1.1}. Lemma \ref
{l2.1} provides bounds for variances of integrals over some exceptional
sets. Lemma \ref{l2.3} gives estimates for variances of integrals over sets
of the form $\left( a,b\right) \cap E_{n}$. Lemma \ref{l2.4} implies bounds
for $\sqrt{n}\int_{E_{n}}\left| {\bf E\,}|f_{\eta }-{\bf E\,}f_{n}|-{\bf E\,}%
|f_{n}-{\bf E\,}f_{n}|\right| $. In Section 3, we replace sets $E_{n}$ by
some sets $C_{n}\subset E_{n}$ removing ''bad'' intervals and tails of small
measure. Then we represent the integral over $C_{n}$ as a sum (in $i$) $%
S_{n} $ of 1-dependent integrals $\delta _{i,n}$ over some sets $I_{i,n}$.
Lemma \ref{l3.3} provides Bernstein-type bounds for moments of summands $%
\delta _{i,n}$. Lemma \ref{l3.4} contains a bound for the correlation
between $S_{n} $ and some centered and normalized Poisson random variable $%
U_{n}=\sum_{i}u_{i,n}$. The summands $u_{i,n}$ are independent centered and
normalized Poisson random variables and the bivariate random vectors $\left(
\delta _{i,n},u_{i,n}\right) $ are 1-dependent. In Lemma \ref{l3.5}, using
bounds from Lemma~\ref{l3.3}, we prove Bernstein-type bounds for moments of
projections of vectors $\left( \delta _{i,n},u_{i,n}\right) $ to
one-dimensional directions. A result of Heinrich \cite{he}, see Lemma \ref
{f3.1}, implies bounds for cumulants of projections of vectors $\left(
S_{n},U_{n}\right) $. In Lemma \ref{l3.6}, we use these bounds to show that
distribution ${\cal L}\left( (S_{n},U_{n})\right) \in {\cal A}_{2}(\tau
_{n}) $ with some $\tau _{n}\leq \tau _{n}^{*}$, where ${\cal A}_{2}(\tau
_{n})$ is a class of distributions introduced by Zaitsev \cite{z86}. In
Section 4, we get bounds for exponential moments of integrals over
exceptional sets, see Lemma \ref{l4.2}. These bounds imply exponential
inequalities for the tails of the corresponding distributions. Theorems \ref
{t1.2}, \ref{t1.3} and \ref{t1.4} are proved in Section 5. We use there a
result of Zaitsev \cite{z02} providing an estimate of the rate of
approximation in a de-Poissonization lemma of Beirlant and Mason \cite{bm}.

\section{\bf Preliminary lemmas.}

\begin{lemma}
[cf. the proof of Gin\'{e}, Mason and Zaitsev \cite{gmz}, Lemma 6.2] \label%
{l2.1}{\it Whenever }$h_{n}\rightarrow 0$ {\it and }$nh_{n}\rightarrow
\infty ${\it , as }$n\rightarrow \infty ,${\it \ for any Borel subset }$B$%
{\it \ of }${\bf R}${\it \ and any sequence of functions }$a_{n}\in L_{1}(%
{\bf R}),$%
\begin{eqnarray}
&&{\bf \,E\,}\left( \sqrt{n}\int_{B}\{|f_{n}(x)-a_{n}(x)|-{\bf E\,}%
|f_{n}(x)-a_{n}(x)|\}\,dx\right) ^{2}  \label{db} \\
&\leq &d(n,B)\stackrel{\rm def}{=}4\,\left\| K\right\| _{\infty }\,{\bf E}\,%
\frac{1}{h_{n}}\int_{B}\left| K\left( \frac{x-X}{h_{n}}\right) \right| \,dx,
\nonumber
\end{eqnarray}
\begin{equation}
{\bf E}\left( \sqrt{n}\int_{B}\{|f_{\eta }(x)-a_{n}(x)|-{\bf E}|f_{\eta
}(x)-a_{n}(x)|\}\,dx\right) ^{2}\leq 2\,d(n,B),  \label{DB}
\end{equation}
{\it where }$d(n,B)$ {\it satisfies}
\begin{equation}
d(n,B)\leq 4\,\kappa ^{2}\,\Omega (n,B)  \label{19a}
\end{equation}
{\it with}
\begin{equation}
\Omega (n,B)\stackrel{\rm def}{=}\left( \int_{B}f(x)\,\,dx+L(n,B)\right) ,
\label{13}
\end{equation}
\begin{equation}
L(n,B)\stackrel{\rm def}{=}\int_{B}\left| h_{n}^{-1}{\bf P}\{X\in
[x-h_{n}/2,x+h_{n}/2]\}-f(x)\right| \,dx\leq L(n,{\bf R})\rightarrow 0,
\label{19b}
\end{equation}
{\it as }${\it \ }n\rightarrow \infty .\medskip $\end{lemma}

\noindent {\it Proof.} Applying the main result in Pinelis \cite{p4}, we get
(see (\ref{k1}))
\begin{eqnarray}
&&{\bf E}\left( \sqrt{n}\int_{B}\{|f_{n}(x)-a_{n}(x)|-{\bf E}%
|f_{n}(x)-a_{n}(x)|\}\,\,dx\right) ^{2}  \nonumber \\
&\leq &4\,{\bf E}\,\left( \frac{1}{h_{n}}\int_{B}\left| K\left( \frac{x-X}{%
h_{n}}\right) \right| \,dx\right) ^{2}  \label{ddss} \\
&\leq &4\,\left\| K\right\| _{\infty }\,{\bf E}\,\frac{1}{h_{n}}%
\int_{B}\left| K\left( \frac{x-X}{h_{n}}\right) \right| \,dx.  \nonumber
\end{eqnarray}
Similarly, taking into account (\ref{kk}) and (\ref{TETTA})--(\ref{Teta}),
we have
\begin{eqnarray*}
&&{\bf E}\left( \sqrt{n}\int_{B}\{|f_{\eta }(x)-a_{n}(x)|-{\bf E}|f_{\eta
}(x)-a_{n}(x)|\}\,\,dx\right) ^{2} \\
&\leq &4\,{\bf E}\left( \frac{1}{h_{n}}\int_{B}\left| \sum_{j\leq \eta
_{1}}K\left( \frac{x-X_{j}}{h_{n}}\right) \right| \,dx\right) ^{2} \\
&\leq &4\,\left\| K\right\| _{\infty }\,{\bf E}\,\frac{1}{h_{n}}%
\int_{B}\left| K\left( \frac{x-X}{h_{n}}\right) \right| \,dx\,{\bf E}\,\eta
_{1}^{2}.
\end{eqnarray*}
Using (\ref{k1}) and (\ref{k2}), we obtain
\[
\left\| K\right\| _{\infty }\,{\bf E}\,\frac{1}{h_{n}}\int_{B}\left| K\left(
\frac{x-X}{h_{n}}\right) \right| \,dx\leq \kappa ^{2}h_{n}^{-1}\int_{B}{\bf P%
}\{X\in [x-h_{n}/2,x+h_{n}/2]\}\,dx.
\]
Furthermore, ${\bf E}\,\eta _{1}^{2}=2$ and
\begin{equation}
h_{n}^{-1}\int_{B}{\bf P}\{X\in [x-h_{n}/2,x+h_{n}/2]\}\,dx\leq
\int_{B}f(x)\,\,dx+L(n,B)=\Omega (n,B).  \label{19x}
\end{equation}
By a special case of Theorem 1 in Chapter 2 of Devroye and Gy\"{o}rfi\ \cite
{dg},
\[
L(n,{\bf R})=\int_{{\bf R}}\left| h_{n}^{-1}\,{\bf P}\{X\in
[x-h_{n}/2,x+h_{n}/2]\}-f(x)\right| \,dx\rightarrow 0,\quad \mbox{as }%
n\rightarrow {\infty },
\]
which completes the proof of Lemma \ref{l2.1}. \medskip

We shall apply Lemma \ref{l2.1} in the case where $a_{n}(x)={\bf E\,}%
f_{n}(x) $. Note that in this situation a similar bound may be derived from
Theorem 2.1 of de Acosta~\cite{a}. Also see Devroye \cite{d}, who obtains
the bound (\ref{ddss}) with $a_{n}(x)=f(x)$. The following standard lemma
follows from Theorem 3 in Chapter 2 of Devroye and Gy\"{o}rfi~\cite{dg}.$%
\medskip $

\begin{lemma}
[see Gin\'{e}, Mason and Zaitsev \cite{gmz},{\bf \ }Lemma 6.1] \label{l2.2}
{\it Suppose that }$H${\it \ is a uniformly bounded real valued} {\it %
function,\ which is equal to zero off a compact interval. Then}
\begin{equation}
|f*H_{h}(x)-I(H)\,f(x)|\rightarrow 0,\quad \mbox{{\it \ as} }h\searrow
0,\quad \mbox{{\it for almost all} }x\in {\bf R},  \label{L2}
\end{equation}
{\it where }$I(H)${\it \ and }$f*H_{h}(x)${\it \ are defined in
}{\rm (\ref {IH})}{\it \ and }{\rm (\ref{9*})}{\it .}$\medskip $
\end{lemma}

\noindent {\it Proof of Lemma }\ref{l1.1}.{\it \ }Applying for each $m\in
{\bf N}$ and for ${\cal H}={\cal H}_{0}$ Lemma 6.1 from Gin\'{e}, Mason and
Zaitsev \cite{gmz}, we conclude that there exist measurable sets $%
Q_{1},Q_{2},\ldots ,Q_{m},\ldots $ such that
\begin{equation}
\int_{Q_{m}}f(x)\,dx\geq 1-2^{-m},  \label{m1}
\end{equation}
$f$\/ is continuous,\ for \ $x\in Q_{m}$, $m=1,2,\ldots ,$ and, uniformly in
$H\in {\cal H}_{0}$,
\begin{equation}
\sup_{x\in Q_{m}}|f*H_{h_{n}}(x)-I(H)\,f(x)|\rightarrow 0,\quad \mbox{as }%
n\rightarrow \infty .  \label{m2}
\end{equation}
Write
\begin{equation}
Q_{s}^{*}=\bigcup_{m=1}^{s}Q_{m}.  \label{m3}
\end{equation}
By (\ref{m1})--(\ref{m3}),
\begin{equation}
\int_{Q_{s}^{*}}f(x)\,dx\geq 1-2^{-s},  \label{m4}
\end{equation}
\begin{equation}
Q_{1}^{*}\subset Q_{2}^{*}\subset \cdots \subset Q_{s}^{*}\subset \cdots
\label{m5}
\end{equation}
and, for $s=1,2,\ldots ,$%
\begin{equation}
\sup_{x\in Q_{s}^{*},\,H\in {\cal H}_{0}}|f*H_{h_{n}}(x)-I(H)\,f(x)|%
\rightarrow 0,\quad \mbox{as }n\rightarrow \infty .  \label{m6}
\end{equation}
Define $m_{1}=1,$%
\begin{equation}
m_{s}=\min \left\{ m>m_{s-1}:\sup_{n\geq m}\;\sup_{x\in Q_{s}^{*},\,H\in
{\cal H}_{0}}\;|f*H_{h_{n}}(x)-I(H)\,f(x)|<2^{-s}\right\} ,  \label{m9}
\end{equation}
for $s=2,3,\ldots $, and
\begin{equation}
F_{l}=Q_{s}^{*},\quad \mbox{for }m_{s}\leq l<m_{s+1}.  \label{m10}
\end{equation}
By (\ref{m4})--(\ref{m10}),
\begin{equation}
\int_{F_{l}}f(x)\,dx\nearrow 1,\quad \mbox{as }l\rightarrow \infty ,
\label{m11}
\end{equation}
\begin{equation}
F_{1}\subset F_{2}\subset \cdots \subset F_{l}\subset \cdots  \label{m12}
\end{equation}
and, for $l=1,2,\ldots ,$%
\begin{equation}
\varepsilon _{l,n}^{*}\stackrel{\rm def}{=}\sup_{m\geq n}\;\sup_{x\in
F_{l},\,H\in {\cal H}_{0}}|f*H_{h_{m}}(x)-I(H)\,f(x)|\rightarrow 0,\quad %
\mbox{as }n\rightarrow \infty .  \label{m13}
\end{equation}

Let sequences $\left\{ \beta _{n}^{*}\right\} _{n=1}^{\infty }$ and $\left\{
D_{n}^{*}\right\} _{n=1}^{\infty }$ satisfy conditions
\begin{equation}
0<\beta _{n}^{*}<D_{n}^{*}<\infty ;\quad \quad \beta _{n}^{*}\searrow
0,\quad D_{n}^{*}\nearrow \infty ,\quad \mbox{as }n\rightarrow \infty .
\label{n3}
\end{equation}
Define, for $l=1,2,\ldots ,$%
\begin{equation}
G_{l}=\left\{ x\in F_{l}:\beta _{l}^{*}\leq f(x)\leq D_{l}^{*}\right\} .
\label{ff}
\end{equation}

Recall that ${\Bbb C}_{n}(x,y)$ and $\rho _{n,x,y}$ were defined in (\ref
{cxy}) and (\ref{ron}). Also observe that
\[
\rho _{n,x,x+th_{n}}=\frac{h_{n}^{-1}{\bf \,E\,}\left[ K\left( \frac{x-X}{%
h_{n}}\right) \,K\left( \frac{x-X}{h_{n}}+t\right) \right] }{\sqrt{h_{n}^{-1}%
{\bf \,E\,}K^{2}\left( \frac{x-X}{h_{n}}\right) h_{n}^{-1}{\bf \,E\,}%
K^{2}\left( \frac{x-X}{h_{n}}+t\right) }},
\]
see (\ref{ron}). Applying Lemma \ref{l2.2}, with $H(u)=K(u)\,K(u+t),$ we
get, for each $t$, that, for almost every $x\in G_{l}$,
\[
h_{n}^{-1}{\bf \,E\,}\left[ K\left( \frac{x-X}{h_{n}}\right) \,K\left( \frac{%
x-X}{h_{n}}+t\right) \right] \rightarrow f(x)\int_{{\bf R}%
}K(u)\,K(u+t)\,du,\quad \mbox{as }n\rightarrow \infty .
\]
Moreover, we get with $H(u)=K^{2}(u)$ and $H(u)=K^{2}(u+t)$, respectively,
for almost every $x\in G_{l},$ both
\[
h_{n}^{-1}{\bf \,E\,}K^{2}\left( \frac{x-X}{h_{n}}\right) \rightarrow
f(x)\,||K^{2}||,\quad \mbox{and\quad }h_{n}^{-1}{\bf \,E\,}K^{2}\left( \frac{%
x-X}{h_{n}}+t\right) \rightarrow f(x)\,||K^{2}||.
\]
Thus, for each $t$ and almost every $x\in G_{l}$,
\[
\rho _{n,x,x+th_{n}}\rightarrow \rho (t),\quad \mbox{as }n\rightarrow \infty
,
\]
and ${\Bbb C}_{n}(x,x+th_{n})\rightarrow {\rm cov}\left( \left| \sqrt{1-\rho
^{2}(t)}Z_{1}+\rho (t)\,Z_{2}\right| ,\left| Z_{2}\right| \right) $. By
Lemma 6.4 from Gin\'{e}, Mason and Zaitsev \cite{gmz}, ${\bf 1}%
_{G_{l}}(x+h_{n}t)\;$converges in measure to$\ {\bf 1}_{G_{l}}(x)=1\;$on $%
G_{l}\times [-1,1]$, and $f(x+h_{n}t)\,{\bf 1}_{G_{l}}(x+h_{n}t)\;$converges
in measure to$\ f(x)\;$on $G_{l}\times [-1,1]$ as functions of $x$ and $t$.
Combining these observations, we readily conclude that $g_{n}(x,t,G_{l})$
converges in measure on $G_{l}\times [-1,1]$ to $g(x,t,G_{l})$. By (\ref{gxt}%
), (\ref{gnxt}), (\ref{kn0}) and (\ref{ff}), functions $g(x,t,G_{l})$ and $%
g_{n}(x,t,G_{l})$ are uniformly bounded on $G_{l}\times [-1,1]$. This
implies
\[
R_{n}(G_{l},G_{l})=\int_{G_{l}}\left( \int_{-1}^{1}\left|
g_{n}(x,t,G_{l})-g(x,t,G_{l})\right| \,dt\right) \,dx\rightarrow 0,\quad %
\mbox{as }n\rightarrow \infty .
\]
It is easy to see that
\begin{equation}
P_{n}\rightarrow 0,\quad \quad \mbox{as }n\rightarrow \infty .  \label{Pn0}
\end{equation}
Define $j_{1}=1$,
\begin{eqnarray}
j_{l} &=&\min \left\{ j>j_{l-1}:\sup_{m\geq j}\left\{ \frac{\sqrt{D_{l}^{*}}%
}{\sqrt{\beta _{l}^{*}}}\left( \frac{1}{\left( \beta _{l}^{*}\,mh_{m}\right)
^{1/5}}+\frac{\varepsilon _{l,m}^{*}}{\beta _{l}^{*}}\right)
+R_{m}(G_{l},G_{l})\right. \right.  \nonumber \\
&&\left. \left. \quad \quad +\;\frac{1}{\beta _{l}^{*}\,\sqrt{mh_{m}^{2}}}%
+\left( \frac{D_{l}^{*}}{\beta _{l}^{*}}\right) ^{3}\,P_{m}<2^{-l}\right\}
\right\} ,\quad \mbox{for }l=2,3,\ldots ,  \label{jl}
\end{eqnarray}
and
\begin{equation}
E_{n}=G_{l}=\left\{ x\in F_{l}:\beta _{l}^{*}\leq f(x)\leq D_{l}^{*}\right\}
,\quad \mbox{for }j_{l}\leq n<j_{l+1}.  \label{EIGL}
\end{equation}
Using (\ref{m11})--(\ref{EIGL}), we obtain
\begin{eqnarray}
\frac{D_{n}^{1/2}}{\beta _{n}^{1/2}}\left( \frac{1}{\left( \beta
_{n}\,nh_{n}\right) ^{1/5}}+\frac{\varepsilon _{n}}{\beta _{n}}\right)
&+&R_{n}(E_{n},E_{n})  \nonumber \\
&+&\frac{1}{\beta _{n}\,\sqrt{nh_{n}^{2}}}+\frac{D_{n}^{3}\,P_{n}}{\beta
_{n}^{3}}\rightarrow 0,\quad \mbox{as }n\rightarrow \infty ,  \label{tendz}
\end{eqnarray}
with
\[
\varepsilon _{n}\leq \sup_{m\geq j_{l}}\varepsilon _{l,m}^{*},\quad \beta
_{n}\geq \beta _{l}^{*},\quad D_{n}\leq D_{l}^{*},\quad \mbox{for }j_{l}\leq
n<j_{l+1}.
\]
It remains to note that, by (\ref{Pc}), (\ref{bdn}) and (\ref{nan}),
\begin{equation}
\beta _{n}\,\lambda (B)\leq {\bf P}(B)\leq D_{n}\,\lambda (B),\quad %
\mbox{for any Borel set }B\subset E_{n},  \label{lamB}
\end{equation}
${\Bbb N}_{n}\leq D_{n}^{1/2}$ and
\begin{equation}
P_{n}\geq c_{f}\,h_{n},  \label{Pn1}
\end{equation}
for sufficiently large $n\geq n_{0}$, where $c_{f}>0$ depends on density $f$
only. Therefore, (\ref{bdn}) and (\ref{tendz}) imply (\ref{tendzero}).
\medskip

\label{r2.1} The choice of the sets $E_{1},E_{2},\ldots ,E_{n},\ldots $
depends on the choice of the sequences $\left\{ \beta _{n}^{*}\right\}
_{n=1}^{\infty }$ and $\left\{ D_{n}^{*}\right\} _{n=1}^{\infty }$ in the
proof of Lemma \ref{l1.1}.$\medskip $ 

In the sequel we shall assume that $h_{n}\rightarrow 0$\ and $%
nh_{n}^{2}\rightarrow \infty $, as $n\rightarrow \infty $ and $n\geq n_{0}$,
where $n_{0}$ is a positive integer which will be chosen as large as it is
necessary for the arguments below to hold. Let $E_{1},E_{2},\ldots
,E_{n},\ldots $ be {\it any }sequence of Borel sets satisfying (\ref{psin0}%
)--(\ref{tendzero}). By (\ref{bdn}) and (\ref{tendzero}), $\frac{\varepsilon
_{n}}{\beta _{n}}\rightarrow 0$ as $n\rightarrow \infty $. Let $n\geq n_{0}$
be so large that
\begin{equation}
\varepsilon _{n}\leq \beta _{n}\,\min \left\{ I(H):H\in {\cal H}_{0}\right\}
/2.  \label{n2}
\end{equation}
Then, by (\ref{bdn}), (\ref{CC}) and (\ref{n2}), for any $x\in E_{n}$, $H\in
{\cal H}_{0}$, we have
\begin{equation}
f(x)\,I(H)/2\leq f*H_{h_{n}}(x)\leq 2\,f(x)\,I(H).  \label{bound}
\end{equation}

We shall use the following fact that follows from Theorem 1 of Sweeting~\cite
{s}.\medskip

\begin{lemma}
\label{f2.1} {\it Let }$(\omega ,\zeta ),(\omega _{1},\zeta _{1}),(\omega
_{2},\zeta _{2}),\ldots ,${\it \ be a sequence of i.i.d. bivariate random
vectors such that each component has variance }$1${\it , mean }$0${\it \ and
finite moments of the third order. Further, let }$(Z_{1}^{*},Z_{2}^{*})${\it %
\ be bivariate normal vector with mean }$0,$ $\mbox{\rm Var}(Z_{1}^{*})=%
\mbox{\rm Var}(Z_{2}^{*})=1,${\it \ and with } $\mbox{\rm cov}%
(Z_{1}^{*},Z_{2}^{*})=\mbox{\rm cov}(\omega ,\zeta )=\rho ${\it . Then there
exists a universal positive constant }$A${\it \ such that }
\begin{equation}
\left| {\bf E\,}\left| \frac{\sum_{i=1}^{n}\zeta _{i}}{\sqrt{n}}\right| -%
{\bf E\,}|Z_{1}^{*}|\right| \leq \frac{A}{\sqrt{n}}\,{\bf E\,}|\zeta |^{3}
\label{in1}
\end{equation}
{\it and, whenever }$\rho ^{2}<1$,
\begin{equation}
\left| {\bf E\,}\left| \frac{\sum_{i=1}^{n}\omega _{i}}{\sqrt{n}}\cdot \frac{%
\sum_{i=1}^{n}\zeta _{i}}{\sqrt{n}}\right| -{\bf E\,}|Z_{1}^{*}Z_{2}^{*}|%
\right| \leq \frac{A}{\left( 1-\rho ^{2}\right) ^{3/2}\sqrt{n}}\left( {\bf %
E\,}|\omega |^{3}+{\bf E\,}|\zeta |^{3}\right)  \label{in2}
\end{equation}
{\it and}
\begin{equation}
\left| {\bf E\,}\left[ \frac{\sum_{i=1}^{n}\omega _{i}}{\sqrt{n}}\cdot
\left| \frac{\sum_{i=1}^{n}\zeta _{i}}{\sqrt{n}}\right| \right] \right| \leq
\frac{A}{\left( 1-\rho ^{2}\right) ^{3/2}\sqrt{n}}\left( {\bf E\,}|\omega
|^{3}+{\bf E\,}|\zeta |^{3}\right) .  \label{in3}
\end{equation}
\medskip \end{lemma}

\begin{lemma}
\label{l2.3} {\it For sufficiently large }$n\geq n_{0}${\it \ and for} {\it %
arbitrary {\rm (}possibly depending on~}$n)$ {\it interval }$\left(
a,b\right) $, $-\infty \leq a<b\leq \infty ,$%
\begin{eqnarray}
&&\left| \,\sigma _{n}^{2}(B)-{\bf P}(B)\,\sigma ^{2}\,\right|  \label{L4} \\
&\leq &A\,{\bf P}(B)\,||K^{2}||\,D_{n}^{1/2}\,\beta _{n}^{-1/2}\,\left(
\frac{\left\| K^{3}\right\| ^{2/5}}{\left\| K^{2}\right\| ^{3/5}\left( \beta
_{n}\,nh_{n}\right) ^{1/5}}+\frac{\varepsilon _{n}}{\left\| K^{2}\right\|
\,\beta _{n}}\right)  \nonumber \\
&&+\;||K^{2}||\,R_{n}(B,E_{n})+16\,\kappa ^{2}\,\left( 1+\beta _{n}^{-1}{\bf %
\,}\varepsilon _{n}\right) \,\min \left\{ P_{n},D_{n}{\bf \,}h_{n}\right\} ,
\nonumber
\end{eqnarray}
{\it where} $B=B\left( n\right) =\left( a,b\right) \cap E_{n}${\it .
Moreover,}
\begin{eqnarray}
&&\left| \,\sigma _{n}^{2}(E_{n})-{\bf P}(E_{n})\,\sigma ^{2}\,\right|
\label{L4a} \\
&\leq &A\,h_{n}^{-1}\,||K^{2}||\,\left( {\Bbb L}_{n}+\frac{\varepsilon _{n}\,%
{\Bbb M}_{n}}{\left\| K^{2}\right\| }\right) +||K^{2}||\,R_{n}(E_{n},E_{n}),
\nonumber
\end{eqnarray}
{\it where\ }${\Bbb L}_{n}${\it \ and }${\Bbb M}_{n}${\it \ are defined in }%
{\rm (\ref{LN})}{\it --}{\rm (\ref{mnab}).}$\medskip $ \end{lemma}

{\it Proof.} Notice that whenever $|x-y|>h_{n}$, random variables $|f_{\eta
}(x)-{\bf E}\,f_{n}(x)|$ and $|f_{\eta }(y)-{\bf E}\,f_{n}(y)|$ are
independent. This follows from the fact that they are functions of
independent increments of the Poisson process with intensity $nf$. Therefore
(see (\ref{TETTA}), (\ref{JnBK}) and (\ref{39a}))
\begin{eqnarray}
v_{n}(B,E_{n}) &=&n\int_{B}\int_{E_{n}}{\bf \,E\,}\{|f_{\eta }(x)-{\bf \,E\,}%
f_{n}(x)|\,|f_{\eta }(y)-{\bf \,E\,}f_{n}(y)|\}\,dx\,dy  \nonumber \\
&&-\;n\int_{B}\int_{E_{n}}\{{\bf \,E\,}|f_{\eta }(x)-{\bf \,E\,}f_{n}(x)|\,%
{\bf E\,}|f_{\eta }(y)-{\bf \,E\,}f_{n}(y)|\}\,dx\,dy  \label{sigm} \\
&=&\int_{B}\int_{E_{n}}{\bf 1}\{|x-y|\leq h_{n}\}\,\mbox{{\rm cov}}\left(
|T_{\eta }(x)|,|T_{\eta }(y)|\right) \sqrt{k_{n}(x)\,k_{n}(y)}\,dx\,dy.
\nonumber
\end{eqnarray}
According to (\ref{sig}) and (\ref{Pc})--(\ref{gnxt}), we have, for $x\in
E_{n}$,
\begin{equation}
\int_{-1}^{1}g(x,t,E_{n})\,dt=\,f(x)\int_{-1}^{1}\mbox{{\rm cov}}\left(
\left| \sqrt{1-\rho ^{2}(t)}Z_{1}+\rho (t)\,Z_{2}\right| ,\left|
Z_{2}\right| \right) \,dt=\frac{f(x)\,\sigma ^{2}}{||K^{2}||},
\label{gxtint}
\end{equation}
\begin{equation}
\int_{B}\left( \int_{-1}^{1}g(x,t,E_{n})\,dt\right) \,dx=\frac{{\bf P}%
(B)\,\sigma ^{2}}{||K^{2}||}  \label{prob}
\end{equation}
and
\begin{equation}
\left| \varphi _{n}^{2}(B)-{\bf P}(B)\,\sigma ^{2}\right| \leq
||K^{2}||\,R_{n}(B,E_{n}),  \label{phips}
\end{equation}
where
\begin{equation}
\varphi _{n}^{2}(B)=||K^{2}||\,\int_{B}\int_{-1}^{1}g_{n}(x,t,E_{n})\,dx\,dt.
\label{fien}
\end{equation}

Furthermore, Var$(Y_{n}(x))=1$ (see (\ref{E}), (\ref{kk}) and (\ref{TETTA}%
)--(\ref{Teta})) and
\begin{equation}
{\bf E\,}|Y_{n}(x)|^{3}\leq A\frac{h_{n}^{-3/2}{\bf \,E\,}\left| K\left(
\frac{x-X}{h_{n}}\right) \right| ^{3}}{\left( h_{n}^{-1}{\bf \,E\,}%
K^{2}\left( \frac{x-X}{h_{n}}\right) \right) ^{3/2}}.  \label{third1}
\end{equation}
Using (\ref{bdn}), (\ref{IH})--(\ref{H0}), (\ref{bound}) and (\ref{third1}),
we get that, for $n\geq n_{0\mbox{,}}$%
\begin{equation}
{\bf E\,}|Y_{n}(x)|^{3}\leq A\frac{2\,\left\| K^{3}\right\| {\bf \,}%
h_{n}^{-1/2}}{\sqrt{f(x)}\,\left( \left\| K^{2}\right\| /2\right) ^{3/2}}%
\leq \frac{A\,\left\| K^{3}\right\| }{\sqrt{\beta _{n}h_{n}}\,\left\|
K^{2}\right\| ^{3/2}}.  \label{third}
\end{equation}

By (\ref{kk}), (\ref{CC}), (\ref{IH}) and (\ref{H0}),
\begin{equation}
\sup_{x\in E_{n}}\left| h_{n}\,k_{n}(x)-\left\| K^{2}\right\| \,f(x)\right|
\leq \varepsilon _{n}.  \label{vbbb}
\end{equation}
Assume that $n\geq n_{0}$ is so large that $\frac{\varepsilon _{n}}{\left\|
K^{2}\right\| \,\beta _{n}}\leq 1/6$, see (\ref{tendzero}).{\bf \ }Thus, for
$x\in E_{n}$, we have
\begin{equation}
h_{n}\,k_{n}(x)=\left\| K^{2}\right\| \,f(x)\,\exp \left( \frac{A\,\theta
\,\varepsilon _{n}}{\left\| K^{2}\right\| \,f(x)}\right) ,  \label{vbb}
\end{equation}
where $|\theta |\leq 1$. Using (\ref{vbb}), we see that, for $x,y\in E_{n}$,
\begin{equation}
\sqrt{k_{n}(x)\,k_{n}(y)}=h_{n}^{-1}\left\| K^{2}\right\| \,\sqrt{f(x)\,f(y)}%
\,\exp \left( \frac{A\,\theta \varepsilon _{n}}{\left\| K^{2}\right\| }%
\left( f^{-1}(x)+f^{-1}(y)\right) \right) .  \label{kn1}
\end{equation}
We shall use the elementary fact that if $X$ and $Y$ are mean zero and
variance $1$ random variables with $\rho ={\bf E}\,XY,$ then $1-{\bf E}%
\,|XY|\leq 1-|\rho |\leq 1-\rho ^{2}$. By an application of Lemma \ref{f2.1}%
, keeping (\ref{Teta}), (\ref{cxy}), (\ref{ron}), (\ref{tendzero}), (\ref{KN}%
), (\ref{KN0}) and (\ref{third}) in mind, we obtain, for $n\geq n_{0}$ large
enough and $x,y\in E_{n}$,
\begin{eqnarray}
&&\left| \,\mbox{{\rm cov}}\left( |T_{\eta }(x)|,|T_{\eta }(y)|\right) -%
{\Bbb C}_{n}(x,y)\right|  \label{kn2} \\
&\leq &A{\bf \,}\min \left\{ 1-\rho _{n,x,y}^{2}+\frac{{\bf E\,}%
|Y_{n}(x)|^{3}+{\bf E\,}|Y_{n}(y)|^{3}}{\sqrt{n}}{\bf \,},\frac{{\bf E\,}%
|Y_{n}(x)|^{3}+{\bf E\,}|Y_{n}(y)|^{3}}{\left( 1-\rho _{n,x,y}^{2}\right)
^{3/2}\sqrt{n}}\right\}  \nonumber \\
&\leq &A{\bf \,}\left( {\Bbb K}_{n}(x,y)+{\Bbb K}_{n}(y,x)\right)  \nonumber
\\
&\leq &\frac{A\,\left\| K^{3}\right\| ^{2/5}\,\left(
f^{-1/5}(x)+f^{-1/5}(y)\right) }{\left\| K^{2}\right\| ^{3/5}\left(
nh_{n}\right) ^{1/5}}\leq \frac{A\,\left\| K^{3}\right\| ^{2/5}}{\left\|
K^{2}\right\| ^{3/5}\left( \beta _{n}\,nh_{n}\right) ^{1/5}}.  \nonumber
\end{eqnarray}
Using (\ref{gnxt}), (\ref{cxy}), (\ref{kn0}), (\ref{bdn}), (\ref{tendzero}),
(\ref{sigm}), (\ref{fien}), (\ref{kn1}), (\ref{kn2}) and the change of
variables $y=x+th_{n}$, we see that, for sufficiently large $n\geq n_{0}$,
\begin{eqnarray}
&&\left| v_{n}(B,E_{n})-\varphi _{n}^{2}(B)\right|  \nonumber \\
&\leq &A\,\int_{B}\int_{E_{n}}{\bf 1}\{|x-y|\leq h_{n}\}\,h_{n}^{-1}\left\|
K^{2}\right\| \,\sqrt{f(x)\,f(y)}\,  \nonumber \\
&&\times \left( {\Bbb K}_{n}(x,y)+{\Bbb K}_{n}(y,x)+\frac{\varepsilon _{n}}{%
\left\| K^{2}\right\| }\,\left( f^{-1}(x)+f^{-1}(y)\right) \right) \,dx\,dy
\label{vphi0} \\
&\leq &A\,{\bf P}(B)\,||K^{2}||\,D_{n}^{1/2}\,\beta _{n}^{-1/2}\,\left(
\frac{\,\left\| K^{3}\right\| ^{2/5}}{\left\| K^{2}\right\| ^{3/5}\left(
\beta _{n}nh_{n}\right) ^{1/5}}+\frac{\varepsilon _{n}}{\left\|
K^{2}\right\| \,\beta _{n}}\right) .  \label{vphi}
\end{eqnarray}

Define
\begin{equation}
\begin{tabular}{ll}
$B_{1}=\left( a-h_{n},a\right) \cap E_{n},\qquad $ & $B_{2}=\left(
b,b+h_{n}\right) \cap E_{n},$ \\
$B_{3}=\left( a,a+h_{n}\right) \cap B,$ & $B_{4}=\left( b-h_{n},b\right)
\cap B.$%
\end{tabular}
\label{BBB}
\end{equation}
Clearly,
\begin{equation}
B=B_{3}\cup B_{4}\cup \left( B\backslash \left( B_{3}\cup B_{4}\right)
\right)  \label{BB}
\end{equation}
and
\begin{equation}
{\bf E\,}J_{n}(B){\bf \,}J_{n}(E_{n}\backslash \left( B\cup B_{1}\cup
B_{2}\right) )=0,  \label{ind0}
\end{equation}
since $J_{n}(B)$ and $J_{n}(E_{n}\backslash \left( B\cup B_{1}\cup
B_{2}\right) )$ are independent. Similarly, according to (\ref{BBB}) and (%
\ref{BB}), ${\bf E\,}J_{n}(B){\bf \,}J_{n}(B_{1}\cup B_{2})={\bf E\,}%
J_{n}(B_{1}){\bf \,}J_{n}(B_{3})+{\bf E\,}J_{n}(B_{2}){\bf \,}J_{n}(B_{4})$.
Note that, by (\ref{CC})--(\ref{H0}), (\ref{19b}) and (\ref{lamB}), we have
\begin{equation}
L(n,B)\leq \lambda (B)\,\varepsilon _{n}\leq \beta _{n}^{-1}\,{\bf P}%
(B)\,\varepsilon _{n},\quad \mbox{for any Borel set }B\subset E_{n}.
\label{LnB}
\end{equation}
By (\ref{JnBK})--(\ref{Pc}), (\ref{PPPP}), (\ref{DB})--(\ref{19b}), (\ref
{lamB}), (\ref{BBB}), (\ref{ind0}) and (\ref{LnB}),
\begin{eqnarray}
\left| \sigma _{n}^{2}(B)-v_{n}(B,E_{n})\right| &=&\left| {\bf E\,}J_{n}(B)%
{\bf \,}J_{n}(B_{1}\cup B_{2})\right|  \nonumber \\
&\leq &\left| {\bf E\,}J_{n}(B_{1}){\bf \,}J_{n}(B_{3})\right| +\left| {\bf %
E\,}J_{n}(B_{2}){\bf \,}J_{n}(B_{4})\right|  \nonumber \\
&\leq &4\,\max_{1\leq i\leq 4}d(n,B_{i})  \label{titt} \\
&\leq &16\,\kappa ^{2}\,\max_{1\leq i\leq 4}\,\left( {\bf P}%
(B_{i})+L(n,B_{i})\right)  \nonumber \\
&\leq &16\,\kappa ^{2}\,\left( 1+\beta _{n}^{-1}{\bf \,}\varepsilon
_{n}\right) \,\max_{1\leq i\leq 4}\,{\bf P}(B_{i})  \nonumber \\
&\leq &16\,\kappa ^{2}\,\left( 1+\beta _{n}^{-1}{\bf \,}\varepsilon
_{n}\right) \,\min \left\{ P_{n},D_{n}{\bf \,}h_{n}\right\} ,  \nonumber
\end{eqnarray}
for sufficiently large $n\geq n_{0}$. Inequalities (\ref{lamB}), (\ref{phips}%
), (\ref{vphi}) and (\ref{titt}) imply (\ref{L4}). Clearly, $\sigma
_{n}^{2}(E_{n})=v_{n}(E_{n},E_{n})$, see (\ref{sigk}). The proof of (\ref
{L4a}) repeats that of (\ref{L4}). Instead of (\ref{vphi}) one should use (%
\ref{vphi0}) coupled with (\ref{mnab}). \medskip

\begin{lemma}
\label{l2.4} {\it For sufficiently large} $n\geq n_{0}$, {\it we have}
\begin{equation}
\int_{E_{n}}\left| {\bf \,}\sqrt{n}{\bf \,E\,}|f_{\eta }(x)-{\bf \,E\,}%
f_{n}(x)|-{\bf E\,}|Z|\sqrt{k_{n}(x)}{\bf \,}\right| {\bf \,}dx\leq \frac{%
A\,\lambda (E_{n})\,\left\| K^{3}\right\| }{\left\| K^{2}\right\| \sqrt{%
nh_{n}^{2}}}  \label{mc1}
\end{equation}
{\it and} $\vspace{-0.2cm}$%
\begin{equation}
\int_{E_{n}}\left| {\bf \,}\sqrt{n}{\bf \,E\,}|f_{n}(x)-{\bf E\,}f_{n}(x)|-%
{\bf E\,}|Z|\sqrt{k_{n}(x)}{\bf \,}\right| \,dx  \label{mc3}
\end{equation}
\[
\leq {\bf \,}\frac{A\,\lambda (E_{n})\,\left\| K^{3}\right\| }{\left\|
K^{2}\right\| \sqrt{nh_{n}^{2}}}+\frac{A\,{\Bbb N}_{n}\sqrt{h_{n}}}{\sqrt{%
\left\| K^{2}\right\| }},
\]
{\it where }${\Bbb N}_{n}${\it \ is defined by }{\rm (\ref{nan}).}\medskip
\end{lemma}

\noindent {\it Proof. }By (\ref{TETTA}), (\ref{Teta}), (\ref{in1}) and (\ref
{third}), for $x\in E_{n}$,
\begin{equation}
\left| \frac{{\bf \,E\,}\left| \sqrt{n}\left\{ f_{\eta }(x)-{\bf E\,}%
f_{n}(x)\right\} \right| }{\sqrt{k_{n}(x)}}-{\bf E\,}|Z|\right| \leq \frac{A%
}{\sqrt{n}}{\bf \,E\,}|Y_{n}(x)|^{3}\leq \frac{A\,\left\| K^{3}\right\| }{%
\sqrt{n\,f(x)\,h_{n}}\,\left\| K^{2}\right\| ^{3/2}}{\bf \,}.  \label{fourth}
\end{equation}
Using (\ref{k3}), (\ref{kk}), (\ref{v}), (\ref{bdn}), (\ref{H0}) and (\ref
{bound}), we get, for $n\geq n_{0}$, $x\in E_{n}$,
\begin{equation}
f(x)\,\left\| K^{2}\right\| \,h_{n}^{-1}/2\leq k_{n}(x)\leq 2\,f(x)\,\left\|
K^{2}\right\| \,h_{n}^{-1}\leq 2\,D_{n}\,\left\| K^{2}\right\| \,h_{n}^{-1}
\label{varb}
\end{equation}
and
\begin{equation}
\left| \sqrt{k_{n}(x)}-\sqrt{n\,\mbox{{\rm Var}}\left( f_{n}(x)\right) }%
\right| \leq \frac{\left( 2\,f(x)\right) ^{2}\sqrt{h_{n}}}{\sqrt{%
f(x)\,\left\| K^{2}\right\| /2}}\leq \frac{A\,f^{3/2}(x)\sqrt{h_{n}}}{\sqrt{%
\left\| K^{2}\right\| }}.  \label{diff}
\end{equation}
Now by (\ref{bdn}), (\ref{tendzero}), (\ref{fourth}) and (\ref{varb}), we
obtain (\ref{mc1}), for sufficiently large $n\geq n_{0}$. Similarly one
obtains
\[
\int_{E_{n}}\left| \,\sqrt{n}{\bf \,E\,}|f_{n}(x)-{\bf E\,}f_{n}(x)|-{\bf E\,%
}|Z|\sqrt{n\,\mbox{{\rm Var}}\left( f_{n}(x)\right) }\,\right| \,dx\leq
\frac{A\,\lambda (E_{n})\,\left\| K^{3}\right\| }{\left\| K^{2}\right\|
\sqrt{nh_{n}^{2}}},
\]
which by (\ref{nan}) and (\ref{diff}) implies (\ref{mc3}). \bigskip

\section{\bf Reduction of the problem to a CLT for 1-dependent random vectors%
}

Let
\begin{equation}
\alpha _{n}\rightarrow 0,\quad \mbox{as }n\rightarrow \infty ,  \label{aln}
\end{equation}
be a non-increasing sequence of strictly positive numbers. In Section 3, we
assume (\ref{aln}) only, keeping in mind that $\alpha _{n}$ will be defined
later by (\ref{aldef}). Using the continuity of our measure, we may find an
interval $[-M_{n},M_{n}]$ so that
\begin{equation}
\alpha _{n}=\int_{\left| x\right| >M_{n}}f(x)\,dx.  \label{n1}
\end{equation}
Assume that $n\geq n_{0}$ is so large that
\begin{equation}
0<\alpha _{n}\leq 1/4\quad \mbox{and}\quad h_{n}\leq \min \left\{
M_{n}/4,1-\alpha _{n}\right\} .  \label{hn2}
\end{equation}
Define $m_{n}=[M_{n}/h_{n}]-1$, $h_{n}^{*}=\left( M_{n}-h_{n}\right) /m_{n}$%
, where $[x]$ denotes the integer part of $x.$ Clearly, by (\ref{hn2}), we
have $M_{n}/2h_{n}\leq m_{n}\leq M_{n}/h_{n}$. Hence,
\begin{equation}
h_{n}\leq h_{n}^{*}\leq 2h_{n}.  \label{hn}
\end{equation}
Recall that $P_{n}$ and $\psi _{n}$ were defined in (\ref{PPPP}) and (\ref
{gamma}). Note that (\ref{tendzero}), (\ref{gamma}), (\ref{aln}) and (\ref
{n1}) imply that ${\bf P}([-M_{n}+h_{n},M_{n}-h_{n}])>\psi _{n}$, for
sufficiently large $n\geq n_{0}$. Define, recurrently, integers $%
l_{1}=-m_{n} $, $l_{i}\in {\bf Z,}$ $l_{1}<l_{2}<\cdot \cdot \cdot
<l_{s_{n}-1}=m_{n}$. Let $l_{i-1}$ be constructed. Then if, for some $l\in
{\bf Z}$, we have ${\bf P}([l_{i-1}h_{n}^{*},\left( l-1\right)
\,h_{n}^{*}])<\psi _{n}$, ${\bf P}([l_{i-1}h_{n}^{*},lh_{n}^{*}])\geq \psi
_{n}$ and ${\bf P}([lh_{n}^{*},M_{n}-h_{n}])\geq \psi _{n}$, we set $l_{i}=l$%
. If, for some $l\in {\bf Z}$, we have ${\bf P}([l_{i-1}h_{n}^{*},\left(
l-1\right) \,h_{n}^{*}])<\psi _{n}$, ${\bf P}([l_{i-1}h_{n}^{*},lh_{n}^{*}])%
\geq \psi _{n}$ and ${\bf P}([lh_{n}^{*},M_{n}-h_{n}])<\psi _{n}$, we set $%
s_{n}-1=i$ and $l_{s_{n}-1}=m_{n}$. Denote
\begin{equation}
z_{0,n}\stackrel{\rm def}{=}-M_{n};\quad \quad z_{s_{n},n}\stackrel{\rm def}{%
=}M_{n};\quad \quad z_{i,n}\stackrel{\rm def}{=}l_{i}h_{n}^{*},\quad \quad %
\mbox{for }i=1,\ldots ,s_{n}-1;  \label{z0n}
\end{equation}
\begin{equation}
I_{i,n}\stackrel{\rm def}{=}E_{n}\cap [z_{i-1,n},z_{i,n}),\quad p_{i,n}%
\stackrel{\rm def}{=}{\bf P}(I_{i,n}),\quad q_{i,n}\stackrel{\rm def}{=}{\bf %
P}([z_{i-1,n},z_{i,n})),  \label{1q}
\end{equation}
for $i=1,\ldots ,s_{n}$. Clearly, we have
\begin{equation}
z_{0,n}<z_{1,n}=-M_{n}+h_{n}<z_{1,n}<\cdot \cdot \cdot
<z_{s_{n}-1,n}=M_{n}-h_{n}<z_{s_{n},n}.  \label{zets}
\end{equation}
Furthermore,
\begin{equation}
P_{n}=\max_{x\in {\bf R}}{\bf P}([x,x+2\,h_{n}])\geq \max_{x\in {\bf R}}{\bf %
P}([x,x+h_{n}^{*}])  \label{Pn}
\end{equation}
(see (\ref{hn})). By (\ref{1q}),
\begin{equation}
p_{i,n}\leq q_{i,n},\quad \quad i=1,\ldots ,s_{n}.  \label{pqin}
\end{equation}
Clearly, by construction, we have
\begin{equation}
\psi _{n}\leq q_{i,n}\leq P_{n}+2\,\psi _{n},\qquad i=2,\ldots ,s_{n}-1,
\label{qin}
\end{equation}
and
\begin{equation}
\max \left\{ q_{1,n},q_{s_{n},n}\right\} \leq P_{n},  \label{qqp}
\end{equation}
for sufficiently large $n\geq n_{0}$. Hence, by (\ref{tendzero}), (\ref
{gamma}), (\ref{Pn0}) and (\ref{pqin})--(\ref{qqp}),
\begin{equation}
\max_{1\leq i\leq s_{n}}\,p_{i,n}\leq \max_{1\leq i\leq
s_{n}}\,q_{i,n}\rightarrow 0,\quad \quad \mbox{as}{\it \ }n\rightarrow
\infty .  \label{777}
\end{equation}

Introduce sets of indices
\begin{equation}
\Upsilon _{1}=\left\{ i=2,\ldots
,s_{n}-1:4\,||K^{2}||\,R_{n}(I_{i,n},E_{n})\geq p_{i,n}\,\sigma ^{2}\right\}
,  \label{Y1}
\end{equation}
\begin{equation}
\Upsilon _{2}=\left\{ i=2,\ldots ,s_{n}-1:p_{i,n}\leq {\bf P}%
([z_{i-1,n},z_{i,n})\backslash I_{i,n})\right\} ,  \label{Y2}
\end{equation}
\begin{equation}
\Upsilon =\Upsilon _{1}\cup \Upsilon _{2},\quad \Upsilon _{3}=\left\{
2,\ldots ,s_{n}-1\right\} \backslash \Upsilon .  \label{Y}
\end{equation}
Define
\begin{equation}
C_{n}=[-M_{n}+h_{n},M_{n}-h_{n}]\cap E_{n}\backslash \bigcup_{i\in \Upsilon
}[z_{i-1,n},z_{i,n}).  \label{cn}
\end{equation}
By construction,
\begin{equation}
C_{n}=\bigcup_{i\in \Upsilon _{3}}I_{i,n},\quad \mbox{and\quad
}I_{i,n}\cap I_{j,n}\hbox{ are empty},\quad \mbox{for }i\neq j.
\label{cunion}
\end{equation}
Using (\ref{gnint}), (\ref{tendzero}), (\ref{1q}), (\ref{zets}) and (\ref{Y1}%
), we obtain
\begin{equation}
{\bf P}\left( \bigcup_{i\in \Upsilon _{1}}I_{i,n}\right) =\sum_{i\in
\Upsilon _{1}}p_{i,n}\leq \frac{4\,||K^{2}||\,R_{n}(E_{n},E_{n})}{\sigma ^{2}%
}\rightarrow 0,\quad \mbox{as }n\rightarrow \infty .  \label{PI1}
\end{equation}
Furthermore, by (\ref{psin0}), (\ref{1q}), (\ref{zets}) and (\ref{Y2}), we
get
\begin{eqnarray}
{\bf P}\left( \bigcup_{i\in \Upsilon _{2}}I_{i,n}\right) &=&\sum_{i\in
\Upsilon _{2}}p_{i,n}\leq \sum_{i\in \Upsilon _{2}}{\bf P}%
([z_{i-1,n},z_{i,n})\backslash I_{i,n})  \nonumber \\
&=&\sum_{i\in \Upsilon _{2}}{\bf P}([z_{i-1,n},z_{i,n})\backslash E_{n})\leq
{\bf P}({\bf R}\backslash E_{n})=\phi _{n}\rightarrow 0,\quad \mbox{as }%
n\rightarrow \infty .  \label{PI2}
\end{eqnarray}
By (\ref{13}), (\ref{19b}), (\ref{Pn0}), (\ref{aln}), (\ref{n1}), (\ref{Pn})
and (\ref{Y1})--(\ref{PI2}), we have
\begin{equation}
\Omega (n,\overline{C}_{n})\leq \alpha _{n}+2\,P_{n}+2\,\phi _{n}+\frac{%
4\,||K^{2}||\,R_{n}(E_{n},E_{n})}{\sigma ^{2}}+L(n,{\bf R)}\rightarrow
0,\quad \mbox{as }n\rightarrow \infty ,  \label{dn0}
\end{equation}
where $\overline{C}_{n}$ denotes the complement of $C_{n}$. By Lemma \ref
{l2.1},
\begin{equation}
{\bf E\,}\left( \sqrt{n}\int_{\overline{C}_{n}}\{|f_{n}(x)-{\bf E\,}%
f_{n}(x)|-{\bf \,E\,}|f_{n}(x)-{\bf E\,}f_{n}(x)|\}\,dx\right) ^{2}\leq d_{n}%
\stackrel{\rm def}{=}d(n,\overline{C}_{n}),  \label{dn1}
\end{equation}
and
\begin{equation}
d_{n}\leq 4\,\kappa ^{2}\,\Omega (n,\overline{C}_{n})\leq 4\,\kappa
^{2}\,\Omega _{n},  \label{dn2}
\end{equation}
where $\Omega _{n}$ is defined in (\ref{dn3}). Similarly, using (\ref{DB})
instead of (\ref{db}), we obtain (see (\ref{JnBK}) and (\ref{sigk}))
\begin{equation}
\sigma _{n}^{2}(E_{n}\backslash C_{n})\leq 8\,\kappa ^{2}\,\left( \alpha
_{n}+2\,P_{n}+\phi _{n}+\frac{4\,||K^{2}||\,R_{n}(E_{n},E_{n})}{\sigma ^{2}}%
+L(n,{\bf R)}\right) \rightarrow 0,  \label{ghhg}
\end{equation}
as $n\rightarrow \infty $. It is easy to see that, by (\ref{bdn}), (\ref{LN}%
), (\ref{mnab}) and (\ref{KN0}),
\begin{equation}
{\Bbb L}_{n}\leq \frac{A\,\left\| K^{3}\right\| ^{2/5}D_{n}^{3/10}\,h_{n}}{%
\left\| K^{2}\right\| ^{3/5}\left( nh_{n}\right) ^{1/5}\,\beta _{n}^{1/2}}%
,\quad {\Bbb M}_{n}\leq 2\,\beta _{n}^{-1}\,h_{n}.  \label{hjk}
\end{equation}
Clearly, ${\bf \,}J_{n}(E_{n})={\bf \,}J_{n}(C_{n})+J_{n}(E_{n}\backslash
C_{n})${\bf . }Therefore, applying (\ref{sigk}), (\ref{psin0}), (\ref{bdn}),
(\ref{tendzero}), (\ref{L4a}), (\ref{cn}), (\ref{ghhg}), (\ref{hjk}) and the
triangle inequality, we get $\sigma _{n}^{2}(C_{n})=\sigma ^{2}+o(1)$ and
\begin{equation}
\frac{1}{2}\,\sigma ^{2}\leq \sigma _{n}^{2}(C_{n})\leq 2\,\sigma ^{2},
\label{p99}
\end{equation}
for sufficiently large $n\geq n_{0}$.

Denote, for $i=1,\ldots ,s_{n}$,
\begin{equation}
\delta _{i,n}\stackrel{\rm def}{=}\frac{\int_{z_{i-1,n}}^{z_{i,n}}{\bf 1}%
_{C_{n}}(x)\,W_{\eta }(x)\,dx}{\sigma _{n}(C_{n})},  \label{2q}
\end{equation}
where
\begin{equation}
W_{\eta }(x)\stackrel{\rm def}{=}\Delta _{\eta }(x)-{\bf E\,}\Delta _{\eta
}(x)=\left( \left| T_{\eta }(x)\right| -{\bf E\,}\left| T_{\eta }(x)\right|
\right) \sqrt{k_{n}(x)},  \label{0q}
\end{equation}
and
\begin{equation}
\Delta _{\eta }(x)\stackrel{\rm def}{=}\sqrt{n}\left| f_{\eta }(x)-{\bf E\,}%
f_{n}(x)\right| =\frac{1}{\sqrt{n}\,h_{n}}\left| \,\sum_{i=1}^{\eta }K\left(
\frac{x-X_{i}}{h_{n}}\right) -n\,{\bf E}\,K\left( \frac{x-X}{h_{n}}\right)
\right| .  \label{deleta}
\end{equation}
Obviously (see (\ref{z0n})--(\ref{zets}), (\ref{Y}) and (\ref{cn})),
\begin{equation}
\delta _{i,n}=0,\quad \mbox{for }i\notin \Upsilon _{3}\qquad \mbox{and}%
\qquad \mbox{ }\delta _{i,n}=\frac{\int_{I_{i,n}}W_{\eta }(x)\,dx}{\sigma
_{n}(C_{n})},\quad \mbox{for }i\in \Upsilon _{3}.  \label{d0}
\end{equation}
Furthermore, $z_{i,n}-z_{i-1,n}\geq h_{n}$, for $i=1,\ldots ,s_{n}$. This
implies that the sequence $\delta _{i,n},$ $1\leq i\leq s_{n}$, is $1$%
-dependent. We used (\ref{k1}), (\ref{deleta}), (\ref{d0}) and that any
functions of the Poisson point process $\left\{ X_{1},\ldots ,X_{\eta
}\right\} $ restricted to disjoint sets are independent.

The use of the sets $C_{n}$ has the advantage over the sets $E_{n}$ in that
they permit us to control the variances of the summands $\delta _{i,n}$ from
below.\medskip

\begin{lemma}
\label{l3.1}{\it \ For sufficiently large }$n>n_{0}$,\ {\it we have}
\[
p_{i,n}\,\sigma ^{2}/4\leq \sigma _{n}^{2}(I_{i,n})\leq 2\,p_{i,n}\,\sigma
^{2},\qquad \mbox{{\it for }}i\in \Upsilon _{3}.
\]
\end{lemma}

\noindent {\it Proof.} According to (\ref{gamma}), (\ref{1q}), (\ref{qin}), (%
\ref{Y2}) and (\ref{Y}), we have, for $i\in \Upsilon _{3}$,
\begin{equation}
p_{i,n}\geq q_{i,n}/2\geq \psi _{n}/2=128\,\kappa ^{2}\,\sigma ^{-2}\,\min
\left\{ P_{n},D_{n}\,h_{n}\right\} .  \label{Ginx}
\end{equation}
Hence, by (\ref{bdn}), (\ref{tendzero}), (\ref{L4}), (\ref{1q}), (\ref{Y1}),
(\ref{Y}) and (\ref{Ginx}), $\beta _{n}^{-1}{\bf \,}\varepsilon _{n}\leq 1$
and
\begin{eqnarray*}
\sigma _{n}^{2}(I_{i,n}) &\geq &p_{i,n}\,\sigma ^{2}-\left| \sigma
_{n}^{2}(I_{i,n})-p_{i,n}\,\sigma ^{2}\right| \\
&\geq &\frac{1}{2}\,p_{i,n}\,\sigma ^{2}-\frac{A\,||K^{2}||\,D_{n}^{1/2}%
\,p_{i,n}}{\beta _{n}^{1/2}}\left( \frac{\left\| K^{3}\right\| ^{2/5}}{%
\left\| K^{2}\right\| ^{3/5}\left( \beta _{n}\,nh_{n}\right) ^{1/5}}+\frac{%
\varepsilon _{n}}{\left\| K^{2}\right\| \,\beta _{n}}\right) \\
&\geq &\frac{1}{4}\,p_{i,n}\,\sigma ^{2},
\end{eqnarray*}
for sufficiently large $n>n_{0}$. Similarly,
\begin{eqnarray*}
\sigma _{n}^{2}(I_{i,n}) &\leq &p_{i,n}\,\sigma ^{2}+\left| \sigma
_{n}^{2}(I_{i,n})-p_{i,n}\,\sigma ^{2}\right| \\
&\leq &\frac{3}{2}\,p_{i,n}\,\sigma ^{2}+\frac{A\,||K^{2}||\,D_{n}^{1/2}%
\,p_{i,n}}{\beta _{n}^{1/2}}\left( \frac{\left\| K^{3}\right\| ^{2/5}}{%
\left\| K^{2}\right\| ^{3/5}\left( \beta _{n}\,nh_{n}\right) ^{1/5}}+\frac{%
\varepsilon _{n}}{\left\| K^{2}\right\| \,\beta _{n}}\right) \\
&\leq &2\,p_{i,n}\,\sigma ^{2},
\end{eqnarray*}
for sufficiently large $n>n_{0}$. \medskip

The following fact will be useful below: if $\xi _{i}$ are independent
centered random variables, then, for every $r\geq 2$,
\begin{equation}
{\bf E\,}\left| \sum_{i=1}^{n}\xi _{i}\right| ^{r}\leq 2^{r+1}\,e^{2r}\,\max
\left[ r^{r/2}\left( \sum_{i=1}^{n}{\bf E\,}\xi _{i}^{2}\right)
^{r/2},\,r^{r}\sum_{i=1}^{n}{\bf E\,}\left| \xi _{i}\right| ^{r}\right]
\label{Pin}
\end{equation}
(Pinelis \cite{p5}, with a unspecified constant $A^{r}$; after
symmetrization, in the form (\ref{Pin}), it follows from Lata\l a \cite{l}).

The following Lemma \ref{l3.2} gives a Rosenthal-type inequality for
Poissonized sums of independent random variables.\medskip

\begin{lemma}
[Gin\'{e}, Mason and Zaitsev \cite{gmz}, Lemma 2.2] \label{l3.2} {\it Assume
that it is known that for any }$n\in {\bf N}$, {\it any i.i.d. centered
random variables }$\xi ,\xi _{1},\xi _{2},\ldots $ {\it \ for some} $r\geq 2$%
,
\begin{equation}
{\bf E\,}\left| \sum_{i=1}^{n}\xi _{i}\right| ^{r}\leq F\left( n\,{\bf E\,}%
\xi ^{2},{\bf \,}n\,{\bf E\,}\left| \xi \right| ^{r}\right) ,  \label{ros}
\end{equation}
{\it where }$F({\bf \,}\cdot \,,{\bf \,}\cdot \,)${\it \ is a non-decreasing
continuous function of two arguments. Then, for any }$\mu >0${\it \ and any
i.i.d. random variables }$\zeta ,\zeta _{1},\zeta _{2},\ldots ,$%
\begin{equation}
{\bf E\,}\left| \sum_{i=1}^{\eta }\zeta _{i}-\mu \,{\bf E\,}\zeta \right|
^{r}\leq F\left( \mu \,{\bf E\,}\zeta ^{2},{\bf \,}\mu {\bf \,E\,}\left|
\zeta \right| ^{r}\right) ,  \label{newros}
\end{equation}
{\it where }$\eta ${\it \ is a Poisson random variable with mean }$\mu ,$%
{\it \ independent of }$\zeta _{1},\zeta _{2},\ldots ${\it .}\medskip
\end{lemma}

\begin{lemma}
\label{l3.3} {\it We have, uniformly in }$i\in \Upsilon _{3}${\it , for
sufficiently large }$n\geq n_{0}${\it \ and for all integers }$r\geq 2$,
\begin{equation}
{\bf E\,}|\delta _{i,n}|^{r}\leq A^{r}\,r^{r}\,p_{i,n}^{r/2-1}\,\left(
\left\| K^{2}\right\| \,D_{n}\,\beta _{n}^{-1}\,\kappa ^{2}\,\sigma
^{-4}\right) ^{r/2}\mbox{\rm Var}(\delta _{i,n}).  \label{D}
\end{equation}
\end{lemma}

\noindent {\it Proof. }By the H\"{o}lder and generalized Minkowski
inequalities (see, e.g., Folland \cite{f}, p. 194), (\ref{0q}) and (\ref{d0}%
),
\begin{equation}
\sigma _{n}^{r}(C_{n})\,{\bf E\,}|\delta _{i,n}|^{r}\leq 2^{r}{\bf E}%
\,\left( \int_{I_{i,n}}\Delta _{\eta }(x)\,dx\right) ^{r}\leq \left(
2\int_{I_{i,n}}\left( {\bf E}\,\Delta _{\eta }^{r}(x)\right)
^{1/r}\,dx\right) ^{r}.  \label{3int}
\end{equation}
Write (see (\ref{deleta}))
\begin{equation}
{\bf E}\,\Delta _{\eta }^{r}(x)=\frac{1}{\left( \sqrt{n}\,h_{n}\right) ^{r}}%
\,{\bf E}\left| {\bf \,}\sum_{i{\bf \,}=1}^{\eta }K\left( \frac{x-X_{i}}{%
h_{n}}\right) -n\,{\bf E}\,K\left( \frac{x-X}{h_{n}}\right) \,\right| ^{r}.
\label{Z1}
\end{equation}
Applying Lemma \ref{l3.2} coupled with inequality (\ref{Pin}), we obtain
\begin{eqnarray}
&&{\bf E}\,\left| {\bf \,}\sum_{i{\bf \,}=1}^{\eta }K\left( \frac{x-X_{i}}{%
h_{n}}\right) -n\,{\bf E}\,K\left( \frac{x-X}{h_{n}}\right) \,\right| ^{r}
\label{Z2} \\
&\leq &2^{r+1}e^{2r}\max \left\{ r^{r/2}\left( n\,{\bf E\,}K^{2}\left( \frac{%
x-X}{h_{n}}\right) \right) ^{r/2},\,r^{r}\,n\,{\bf E}\,\left| K\left( \frac{%
x-X}{h_{n}}\right) \right| ^{r}\right\} .  \nonumber
\end{eqnarray}
Therefore, using (\ref{k2}), (\ref{IH})--(\ref{H0}), (\ref{bound}), (\ref{cn}%
), (\ref{Z1}) and (\ref{Z2}), we see that, for $n\geq n_{0}$, $x\in C_{n}$,
the moment ${\bf E}\,\Delta _{\eta }^{r}(x)$ may be estimated from above by
\begin{eqnarray*}
&&\frac{2^{r+1}\,e^{2r}}{\left( \sqrt{n}\,h_{n}\right) ^{r}}\,\max \left\{
r^{r/2}\left( n\,{\bf E\,}K^{2}\left( \frac{x-X}{h_{n}}\right) \right)
^{r/2},\,r^{r}\,n\,{\bf E\,}\left| K\left( \frac{x-X}{h_{n}}\right) \right|
^{r}\right\} \\
&\leq &\!2^{r+1}e^{2r}\,\max \left\{ r^{r/2}\left( 2\,f(x)\,\left\|
K^{2}\right\| \,h_{n}^{-1}\right) ^{r/2}\!\!,2\,r^{r}n^{1-r/2}\kappa
^{r-2}f(x)\left\| K^{2}\right\| h_{n}^{1-r}\right\} \!.
\end{eqnarray*}
Since $\beta _{n}nh_{n}\rightarrow \infty $, as $n\rightarrow \infty $ (see (%
\ref{bdn}) and (\ref{tendzero})), we estimate for sufficiently large $n\geq
n_{0}$, $x\in C_{n}$,
\[
{\bf E}\,\Delta _{\eta }^{r}(x)\leq 2^{r+1}\,e^{2r}\,r^{r}\left(
2\,f(x)\,\left\| K^{2}\right\| \,h_{n}^{-1}\right) ^{r/2}.
\]
Substituting this into (\ref{3int}), and using H\"{o}lder's inequality, we
get
\begin{equation}
{\bf E\,}|\delta _{i,n}|^{r}\leq A^{r}\,r^{r}\sigma
_{n}^{-r}(C_{n})\,\lambda ^{r/2}(I_{i,n})\,\left( p_{i,n}\,\left\|
K^{2}\right\| \,h_{n}^{-1}\right) ^{r/2},  \label{momm}
\end{equation}
where $p_{i,n}$ is defined in (\ref{1q}). By Lemma \ref{l3.1},
\begin{equation}
\sigma _{n}^{2}(I_{i,n})\geq p_{i,n}\,\sigma ^{2}/4,  \label{sigG}
\end{equation}
for sufficiently large $n\geq n_{0}$. It is easy to see that
\begin{equation}
\mbox{Var}(\delta _{i,n})=\frac{\sigma _{n}^{2}(I_{i,n})}{\sigma
_{n}^{2}(C_{n})}.  \label{vard}
\end{equation}
Each $I_{i,n}$, $i=2,\ldots ,s_{n}-1$, can be represented as $I_{i,n}=\left(
J_{i,n}\cup L_{i,n}\right) \cap E_{n}$, where $J_{i,n}$ is an interval of
length $h_{n}^{*}$ and $L_{i,n}$ is a set with ${\bf P}(L_{i,n})\leq 2\,\psi
_{n}$ with $\psi _{n}$ defined in (\ref{gamma}). Therefore, by (\ref{sikap}%
), (\ref{bdn}), (\ref{gamma}), (\ref{lamB}), (\ref{hn}) and (\ref{1q}),
\begin{eqnarray}
\lambda (I_{i,n}) &\leq &\lambda (J_{i,n}\cap E_{n})+\lambda (L_{i,n}\cap
E_{n})\leq 2\,h_{n}+\beta _{n}^{-1}\,{\bf P}(L_{i,n})  \nonumber \\
&\leq &\left( 2+512\,D_{n}\,\beta _{n}^{-1}\,\kappa ^{2}\,\sigma
^{-2}\right) \,h_{n}\leq A\,D_{n}\,\beta _{n}^{-1}\,\kappa ^{2}\,\sigma
^{-2}h_{n}.  \label{uuuhh}
\end{eqnarray}
Substituting (\ref{sigG}) into (\ref{momm}) and using (\ref{p99}), (\ref
{vard}) and (\ref{uuuhh}), we obtain inequality (\ref{D}). \medskip

Define
\begin{equation}
S_{n}=\sum_{i=1}^{s_{n}}\delta _{i,n}=\sum_{i\in \Upsilon _{3}}\delta _{i,n}=%
\frac{\int_{C_{n}}W_{\eta }(x)\,dx}{\sigma _{n}(C_{n})}  \label{A}
\end{equation}
(see (\ref{Y}), (\ref{cn}) and (\ref{d0})),
\begin{equation}
U_{n}=\frac{1}{\sqrt{n}}\left\{ \sum_{j\leq \eta }{\bf 1}\{X_{j}\in \left[
-M_{n},M_{n}\right] \}-n{\bf \,P}\{X\in \left[ -M_{n},M_{n}\right] \}\right\}
\label{u}
\end{equation}
and
\begin{equation}
V_{n}=\frac{1}{\sqrt{n}}\left\{ \sum_{j\leq \eta }{\bf 1}\{X_{j}\notin
\left[ -M_{n},M_{n}\right] \}-n{\bf \,P}\{X\notin \left[ -M_{n},M_{n}\right]
\}\right\} .  \label{Vn}
\end{equation}
Set
\begin{equation}
u_{i,n}=\frac{1}{\sqrt{n}}\left\{ \sum_{j\leq \eta }{\bf 1}\left\{ X_{j}\in
[z_{i-1,n},z_{i,n})\right\} -n\,q_{i,n}\right\} ,\quad \quad i=1,\ldots
,s_{n}.  \label{3q}
\end{equation}
It is easy to see that $\sqrt{n}\,u_{i,n}$ is a centered Poisson random
variable with
\begin{equation}
\mbox{Var}(\sqrt{n}\,u_{i,n})=n\,q_{i,n},\quad \quad i=1,\ldots ,s_{n}.
\label{varu}
\end{equation}
Recall that we have $C_{n}\subset \left[ -M_{n}+h_{n},M_{n}-h_{n}\right] $,
see (\ref{cn}). Clearly, $(S_{n},U_{n})$ is a function of the Poisson point
process $\left\{ X_{1},\ldots ,X_{\eta }\right\} $ restricted to the set $%
\left[ -M_{n},M_{n}\right] $ and $V_{n}$ is a function of the same process
restricted to the set ${\bf R}\backslash \left[ -M_{n},M_{n}\right] $.
Therefore, $(S_{n},U_{n})$ is independent of $V_{n}$. Obviously,
\[
U_{n}=\sum_{i=1}^{s_{n}}u_{i,n}
\]
and summands $u_{i,n}$, $i=1,\ldots ,s_{n}$, are independent. Hence,
\begin{equation}
\mbox{{\rm Var}}\left( U_{n}\right) =\sum_{i=1}^{s_{n}}\mbox{{\rm Var}}%
\left( u_{i,n}\right) =\sum_{i=1}^{s_{n}}q_{i,n}={\bf P}\{X\in \left[
-M_{n},M_{n}\right] \},  \label{varuu}
\end{equation}
see (\ref{z0n})--(\ref{zets}). Observe that
\begin{equation}
\mbox{{\rm Var}}(S_{n})=1\quad \mbox{and}\quad \mbox{{\rm Var}}%
(U_{n})=1-\alpha _{n},  \label{vars}
\end{equation}
where $\alpha _{n}={\bf P}\{X\notin \left[ -M_{n},M_{n}\right] \}$, see (\ref
{n1})$.$ \medskip

\begin{lemma}
\label{l3.4} {\it For sufficiently large} $n\geq n_{0}$, {\it we have} $%
\ldots $%
\begin{equation}
\left| \,\mbox{{\rm cov}}(S_{n},U_{n})\right| \leq \frac{A\,\left\|
K^{3}\right\| \,\lambda (E_{n})}{\sigma \,\left\| K^{2}\right\| \sqrt{%
nh_{n}^{2}}}.  \label{cov}
\end{equation}
{\it Moreover,}
\begin{equation}
\max_{i\in \Upsilon _{3}}\frac{\left| \,\mbox{{\rm cov}}(\delta
_{i,n},u_{i,n})\right| }{\left( \mbox{{\rm Var}}(u_{i,n})\mbox{{\rm Var}}%
(\delta _{i,n})\right) ^{1/2}}\rightarrow 0,{\it \quad as\ }n\rightarrow
\infty .  \label{corr}
\end{equation}
\end{lemma}
\medskip

\noindent {\it Proof.} According to (\ref{1q}), (\ref{0q}), (\ref{d0}) and (%
\ref{3q}), we have, for $i\in \Upsilon _{3}$,
\begin{equation}
\sigma _{n}(C_{n})\,\mbox{{\rm cov}}(\delta
_{i,n},u_{i,n})=q_{i,n}^{1/2}\int_{\,I_{i,n}}\left( {\bf E\,}\left| T_{\eta
}(x)\right| \,u_{i,n}\,q_{i,n}^{-1/2}\right) \sqrt{k_{n}(x)}\,dx.
\label{bubu}
\end{equation}
Note that (\ref{qin}) and (\ref{777}) imply that
\begin{equation}
\psi _{n}\leq \min_{2\leq i\leq s_{n}-1}q_{i,n}\rightarrow 0,{\it \quad }%
\mbox{as}{\it \ }n\rightarrow \infty .  \label{qin0}
\end{equation}
Below we assume that $n\geq n_{0}$ is sufficiently large. By (\ref{lamB}), (%
\ref{1q}), (\ref{sigG}) and (\ref{vard}),
\begin{equation}
\lambda (I_{i,n})\leq \beta _{n}^{-1}\,p_{i,n}\leq 4\,\beta
_{n}^{-1}\,\sigma ^{-2}\sigma _{n}^{2}(I_{i,n})=4\,\beta _{n}^{-1}\,\sigma
^{-2}\sigma _{n}^{2}(C_{n})\,\mbox{{\rm Var}}(\delta _{i,n}).  \label{bububu}
\end{equation}
Note now that
\begin{equation}
\left( T_{\eta }(x),u_{i,n}\,q_{i,n}^{-1/2}\right)
=_{d}n^{-1/2}\sum_{l=1}^{n}\left( Y_{n}^{(l)}(x),U^{(l)}\right) ,
\label{354a}
\end{equation}
where $\left( Y_{n}^{(l)}(x),U^{(l)}\right) ,$ $l=1,\ldots ,n$, are i.i.d. $%
\left( Y_{n}(x),U\right) ,$ with $Y_{n}(x)$ defined in (\ref{Yn}) and
\begin{equation}
U=q_{i,n}^{-1/2}\left\{ \sum_{j\leq \eta _{1}}{\bf 1}\left\{ X_{j}\in
[z_{i-1,n},z_{i,n})\right\} -q_{i,n}\right\} ,  \label{zaq}
\end{equation}
$\eta _{1}$ denoting a Poisson random variable with mean $1$ from (\ref{Yn}%
), which is independent of $X,X_{1},X_{2},\dots .$ Using (\ref{k1}), (\ref
{sikap}), (\ref{Yn}), (\ref{IH})--(\ref{H0}), (\ref{gamma}), (\ref{bound}), (%
\ref{qin}) and (\ref{zaq}), we see that, for any $x\in C_{n}$,
\begin{equation}
\left| \mbox{{\rm cov}}(Y_{n}(x),U)\right| =\frac{\left| {\bf E\,}\left[
K\left( \frac{x-X}{h_{n}}\right) \,{\bf 1}\left\{ X\in
[z_{i-1,n},z_{i,n})\right\} \right] \right| }{q_{i,n}^{1/2}{\bf \,}\left(
{\bf E\,}K^{2}\left( \frac{x-X}{h_{n}}\right) \right) ^{1/2}}\leq \frac{2%
\sqrt{2\,D_{n}}\,\kappa \,h_{n}^{1/2}}{q_{i,n}^{1/2}{\bf \,}\left\|
K^{2}\right\| ^{1/2}}\leq \frac{1}{4},  \label{correl}
\end{equation}
if $\psi _{n}=256\,\kappa ^{2}\,\sigma ^{-2}\,D_{n}\,h_{n}$. Furthermore,
using the first equality in (\ref{correl}), (\ref{k1}), (\ref{sikap}), (\ref
{PPPP}), (\ref{qin}) and H\"{o}lder's inequality, we get
\begin{equation}
\left| \mbox{{\rm cov}}(Y_{n}(x),U)\right| \leq q_{i,n}^{-1/2}{\bf \,P}%
^{1/2}\left\{ X\in [x-h_{n}/2,x+h_{n}/2]\right\} \leq \psi _{n}^{-1/2}{\bf \,%
}P_{n}^{1/2}\leq \frac{1}{8\sqrt{2}},  \label{correl2}
\end{equation}
if $\psi _{n}=256\,\kappa ^{2}\,\sigma ^{-2}\,P_{n}$.

Applying part (\ref{in3}) of Lemma \ref{f2.1} and using (\ref{third}), (\ref
{354a}), (\ref{correl}), (\ref{correl2}) and inequality (\ref{newros}) of
Lemma \ref{l3.2} in the case ${\bf P}\left\{ \zeta =1\right\} =1$ together
with inequality (\ref{Pin}), we get
\begin{eqnarray}
{\bf E\,}\left[ \left| T_{\eta }(x)\right| \,u_{i,n}\,q_{i,n}^{-1/2}\right]
&\leq &\frac{A}{\sqrt{n}}\,\left( {\bf E\,}\left| Y_{n}(x)\right| ^{3}+{\bf %
E\,}\left| U\right| ^{3}\right)  \nonumber \\
&\leq &\frac{A}{\sqrt{n}}\,\left( \frac{\left\| K^{3}\right\| }{\left\|
K^{2}\right\| ^{3/2}\sqrt{f(x)\,h_{n}}}+\frac{1}{\sqrt{q_{i,n}}}\right) .
\label{bu}
\end{eqnarray}
Using (\ref{bdn}), (\ref{tendzero}), (\ref{gamma}), (\ref{Pn1}), (\ref{varb}%
), (\ref{p99}), (\ref{uuuhh}), (\ref{varu}), (\ref{bubu})--(\ref{bububu})
and (\ref{bu}), we get (\ref{corr}):
\begin{eqnarray*}
&&\max_{i\in \Upsilon _{3}}\frac{\left| \,\mbox{{\rm cov}}(\delta
_{i,n},u_{i,n})\right| }{\left( \mbox{{\rm Var}}(u_{i,n})\,\mbox{{\rm Var}}%
(\delta _{i,n})\right) ^{1/2}} \\
&\leq &A\,\max_{i\in \Upsilon _{3}}\,\frac{q_{i,n}^{1/2}\sqrt{\lambda \left(
I_{i,n}\right) \,\lambda \left( I_{i,n}\right) }}{\sigma _{n}(C_{n})\left(
q_{i,n}\,\mbox{{\rm Var}}(\delta _{i,n})\right) ^{1/2}} \\
&&\times \max_{x\in E_{n}}\left\{ \left( \frac{\left\| K^{3}\right\| }{%
\left\| K^{2}\right\| ^{3/2}\sqrt{f(x)\,h_{n}}}+\frac{1}{\sqrt{q_{i,n}}}%
\right) \sqrt{\frac{f(x)\left\| K^{2}\right\| }{n\,h_{n}}}\right\} \\
&\leq &\frac{A\,\left( D_{n}\,\beta _{n}^{-1}\,\kappa ^{2}\,\sigma
^{-2}\right) ^{1/2}}{\sigma \,\beta _{n}^{1/2}\sqrt{n}} \\
&&\times \max_{x\in E_{n}}\left\{ \left( \frac{\left\| K^{3}\right\| }{%
\left\| K^{2}\right\| ^{3/2}\sqrt{f(x)\,h_{n}}}+\frac{1}{\sqrt{\psi _{n}}}%
\right) \sqrt{f(x)\left\| K^{2}\right\| }\right\} \rightarrow
0,\quad \hbox{ as }n\rightarrow \infty .
\end{eqnarray*}
Similarly,
\begin{equation}
\sigma _{n}(C_{n})\,\,\mbox{{\rm cov}}(S_{n},U_{n})=(1-\alpha
_{n})^{1/2}\int_{C_{n}}{\bf E\,}\left[ \left| T_{\eta }(x)\right|
\,U_{n}\,(1-\alpha _{n})^{-1/2}\right] \sqrt{k_{n}(x)}\,dx  \label{nuh}
\end{equation}
and
\[
\left| \,\mbox{{\rm cov}}(T_{\eta }(x),U_{n}\,(1-\alpha
_{n})^{-1/2})\,\right| \leq 1/4.
\]
Applying part (\ref{in3}) of Lemma \ref{f2.1} and using (\ref{Teta}), (\ref
{third}) and again inequality (\ref{newros}) of Lemma \ref{l3.2} in the case
${\bf P}\left\{ \zeta =1\right\} =1$ coupled with inequality (\ref{Pin}), we
get
\begin{equation}
{\bf E\,}\left[ \left| T_{\eta }(x)\right| \,U_{n}\,(1-\alpha
_{n})^{-1/2}\right] \leq \frac{A}{\sqrt{n}}\,\left( \frac{\left\|
K^{3}\right\| }{\left\| K^{2}\right\| ^{3/2}\sqrt{f(x)\,h_{n}}}+\frac{1}{%
\sqrt{1-\alpha _{n}}}\right) .  \label{TU}
\end{equation}
By (\ref{cn}),
\begin{equation}
\lambda (C_{n})\leq \lambda (E_{n}).  \label{nunu}
\end{equation}
Using (\ref{bdn}), (\ref{tendzero}), (\ref{varb}), (\ref{hn2}), (\ref{p99}),
(\ref{nuh})--(\ref{nunu}), we get (\ref{cov}):
\begin{eqnarray*}
&&\left| \,\mbox{{\rm cov}}(S_{n},U_{n})\right| \leq \frac{A\,(1-\alpha
_{n})^{1/2}\,\lambda (E_{n})}{\sigma _{n}(C_{n})}{\bf \,} \\
&&\times \;\frac{1}{\sqrt{n}}\,\max_{x\in E_{n}}\left\{ \left( \frac{\left\|
K^{3}\right\| }{\left\| K^{2}\right\| ^{3/2}\sqrt{f(x)\,h_{n}}}+\frac{1}{%
\sqrt{1-\alpha _{n}}}\right) \sqrt{f(x)\,\left\| K^{2}\right\| \,h_{n}^{-1}}%
\right\} \\
&\leq &\frac{A\,\left\| K^{3}\right\| \,\lambda (E_{n})}{\sigma \,\left\|
K^{2}\right\| \sqrt{nh_{n}^{2}}}.
\end{eqnarray*}
\medskip

Below, for $z=\left( z_{1},z_{2}\right) $, $u=\left( u_{1},u_{2}\right) \in
{\bf C}^{2}$, we shall use the notation
\[
\left| z\right| =\left| z_{1}\right| +\left| z_{2}\right| ,\quad \left\|
z\right\| ^{2}=\left| z_{1}\right| ^{2}+\left| z_{2}\right| ^{2},\quad
\left\langle z,u\right\rangle =z_{1}\overline{u_{1}}+z_{2}\overline{u_{2}}.
\]
We shall write $\Gamma _{r}\left\{ \xi \right\} $ for the $k$-th cumulant of
a random variable $\xi $. Recall that if, for some $c>0$, a random variable $%
\xi $ has finite exponential moments ${\bf E\,}e^{z\xi }$, $z\in {\bf C}$, $%
|z|<c$, then (choosing $\log 1=0$)
\begin{equation}
\log {\bf E\,}e^{z\xi }=\sum_{r=0}^{\infty }\frac{\Gamma _{r}\left\{ \xi
\right\} \,z^{r}}{r!}\quad \mbox{and\quad }\Gamma _{r}\left\{ \xi \right\}
=\left. \frac{d^{r}}{dz^{r}}\log {\bf E\,}e^{z\xi }\right| _{z=0}.
\label{cum1}
\end{equation}
Clearly, $\Gamma _{0}\left\{ \xi \right\} =0$, $\Gamma _{1}\left\{ \xi
\right\} ={\bf E\,}\xi $, $\Gamma _{2}\left\{ \xi \right\} =$Var$\left( \xi
\right) ,$%
\begin{equation}
\Gamma _{r}\left\{ a\xi \right\} =a^{r}\Gamma _{r}\left\{ \xi \right\}
,\quad r=0,1,\ldots .  \label{aG}
\end{equation}
In the two-dimensional case, when $\xi =\left( \xi _{1},\xi _{2}\right) $ is
a bivariate random vector, if $\left| {\bf E\,}e^{\left\langle z,\xi
\right\rangle }\right| <\infty $, $z\in {\bf C}^{2}$, $|z|<c$, $c>0$, then
\begin{equation}
\log {\bf E\,}e^{\left\langle z,\xi \right\rangle
}=\!\!\!\sum_{r_{1},r_{2}=0}^{\infty }\!\!\!\!\frac{\Gamma
_{r_{1},r_{2}}\left\{ \xi \right\} \,z_{1}^{r_{1}}\,z_{2}^{r_{2}}}{%
r_{1}!\,r_{2}!},\ \mbox{where }\Gamma _{r_{1},r_{2}}\left\{ \xi \right\}
=\left. \frac{\partial ^{r_{1}+r_{2}}}{\partial z_{1}^{r_{1}}\,\partial
z_{2}^{r_{2}}}\log {\bf E\,}e^{\left\langle z,\xi \right\rangle }\right|
_{z=0}.  \label{cum2}
\end{equation}
\medskip

\begin{lemma}
[a particular case of Heinrich \cite{he}, Lemma 5] \label{f3.1} {\it Let }$%
\zeta _{1},\zeta _{2},\ldots ,\zeta _{m}$ {\it be }$1${\it -dependent
bivariate random vectors with zero means. Let }$\Lambda _{i}^{2}${\it \ be
the maximal eigenvalue of the covariance matrix of }$\zeta _{i}${\it , }$%
i=1,\ldots ,m${\it . Let }$\lambda ^{2}${\it \ be the minimal eigenvalue of
the covariance matrix }${\bf B}${\it \ of }$\Xi =\zeta _{1}+\zeta _{2}+\cdot
\cdot \cdot +\zeta _{m}${\it . Set }$\Theta ={\bf B}^{-1/2}\Xi ${\it .
Assume that there exists a constant }$H\geq 1/2${\it \ and a real number }$%
\gamma ${\it \ such that}
\begin{equation}
18\,H\,\max_{1\leq i\leq m}\Lambda _{i}^{2}\leq \gamma ^{2}  \label{HL}
\end{equation}
{\it and, for any }$t\in {\bf R}^{2},$%
\begin{equation}
\left| {\bf E\,}\left\langle t,\zeta _{i}\right\rangle ^{r}\right| \leq
H\,r!\,\gamma ^{r-2}\,\left| t\right| ^{r-2}\,\mbox{\rm Var}\left(
\left\langle t,\zeta _{i}\right\rangle \right) ,\quad i=1,\ldots ,m,\quad
r=3,4,\ldots .  \label{BH}
\end{equation}
{\it Then}
\begin{equation}
\sup_{\left\| t\right\| =1}\left| \Gamma _{r}\left\{ \left\langle t,\Theta
\right\rangle \right\} \right| \leq H^{*}\left( r-2\right) !\left( 8\sqrt{2}%
\,\gamma /\lambda \right) ^{r-2},\quad r=2,3,\ldots ,  \label{cum}
\end{equation}
{\it where }$H^{*}=280\,H\,\lambda ^{-2}\sum_{i=1}^{m}\Lambda
_{i}^{2}.\medskip $ \end{lemma}

Note that (\ref{BH}) is automatically satisfied for $r=2$, since $H\geq 1/2$%
.\medskip

\begin{lemma}
\label{l3.5} {\it For sufficiently large }$n\geq n_{0}$, {\it we have,} {\it %
uniformly in }$i=1,\ldots ,s_{n}${\it ,}
\begin{equation}
{\bf E\,}\left| t_{1\,}\delta _{i,n}+t_{2}\,u_{i,n}\right| ^{r}\leq
A\,r!\,\gamma _{n}^{r-2}\,\left\| t\right\| ^{r-2}\,\mbox{\rm Var}%
(t_{1}\,\delta _{i,n}+t_{2}\,u_{i,n}),  \label{Bern}
\end{equation}
{\it for all integers }$r\geq 2$ {\it and for all }$t=\left(
t_{1},t_{2}\right) \in {\bf R}^{2}$, {\it where}
\begin{equation}
\gamma _{n}=A\,\left( \Psi _{n}^{3/2}\,\,\max_{i\in \Upsilon
_{3}}\,p_{i,n}^{1/2}+\max_{1\leq i\leq s_{n}}\,q_{i,n}^{1/2}\right)
\rightarrow 0{\it ,\ as\ }n\rightarrow \infty ,  \label{fla}
\end{equation}
{\it and }$\Psi _{n}${\it \ is defined in }{\rm (\ref{PSI})}. {\it Moreover,
for all integers }$r\geq 3$,
\begin{equation}
\sup_{\left\| t\right\| =1}\left| \Gamma _{r}\left\{
t_{1}\,S_{n}+t_{2}\,U_{n}\right\} \right| \leq \left( r-2\right) !\left(
A\,\gamma _{n}\right) ^{r-2}.  \label{hei8}
\end{equation}
\end{lemma}

\noindent {\it Proof. }Let us prove (\ref{Bern}). Without loss of generality
we assume that
\begin{equation}
\left\| t\right\| =1.  \label{t1t2}
\end{equation}
Applying inequality (\ref{newros}) of Lemma \ref{l3.2} in the case ${\bf P}%
\left\{ \zeta =1\right\} =1-{\bf P}\left\{ \zeta =0\right\} =q_{i,n}$ (see (%
\ref{3q})) coupled with inequality (\ref{Pin}), we get, for $i=1,\ldots
,s_{n}$,
\begin{equation}
{\bf E\,}\left| u_{i,n}\right| ^{r}\leq A^{r}\,n^{-r/2}\,\left(
r^{r/2}\,\left( n\,q_{i,n}\right) ^{r/2}+r^{r}\,n\,q_{i,n}\right) .
\label{Eur}
\end{equation}
Using (\ref{varu}) and (\ref{Eur}), we obtain
\begin{equation}
{\bf E\,}\left| u_{i,n}\right| ^{r}\leq A^{r}\,r^{r}\,\left(
q_{i,n}+n^{-1}\right) ^{r/2-1}\,\mbox{Var}(u_{i,n}).  \label{eur}
\end{equation}
Relation (\ref{corr}) of Lemma \ref{l3.4} implies that
\begin{eqnarray}
\mbox{Var}(t_{1}\,\delta _{i,n}+t_{2}\,u_{i,n}) &=&t_{1}^{2}\,\mbox{Var}%
(\delta _{i,n})+t_{2}^{2}\,\mbox{Var}(u_{i,n})+2\,t_{1}t_{2}\,%
\mbox{{\rm
cov}}(\delta _{i,n},u_{i,n})  \nonumber \\
&\geq &\frac{1}{2}\,\left( t_{1}^{2}\,\mbox{Var}(\delta _{i,n})+t_{2}^{2}\,%
\mbox{Var}(u_{i,n})\right) ,  \label{eur1}
\end{eqnarray}
if $n\geq n_{0}$ is large enough (for $i\notin \Upsilon _{3}$ inequality (%
\ref{eur1}) is trivial, see (\ref{d0})). Recall that $nh_{n}^{2}\rightarrow
\infty $, as $n\rightarrow \infty $. Therefore, (\ref{gamma}), (\ref{Pn1})
and (\ref{qin}) imply that
\begin{equation}
n^{-1}\leq q_{i,n},\quad \mbox{for }i=2,\ldots ,s_{n}-1  \label{nq}
\end{equation}
and sufficiently large $n\geq n_{0}$. Notice that $y\leq \left( y+1\right)
^{r-2},$ for $y\geq 0,$ $r\geq 2$. Moreover, by (\ref{sikap}), (\ref{bdn})
and (\ref{PSI}), we have $\Psi _{n}\geq 1/4$. Hence, applying Lemma~\ref
{l3.3} together with (\ref{PSI}), (\ref{d0}), (\ref{fla}) and (\ref{t1t2})--(%
\ref{eur1}), we get (\ref{Bern}):
\begin{eqnarray}
&&{\bf E\,}\left| t_{1\,}\delta _{i,n}+t_{2}\,u_{i,n}\right| ^{r}
\label{gmm} \\
&\leq &2^{r}\,{\bf E\,}\left| t_{1}\,\delta _{i,n}\right| ^{r}+2^{r}\,{\bf %
E\,}\left| t_{2}\,u_{i,n}\right| ^{r}  \nonumber \\
&\leq &A^{r}\,r^{r}\left( p_{i,n}^{r/2-1}\,\left( \left\| K^{2}\right\|
\,D_{n}\,\beta _{n}^{-1}\,\kappa ^{2}\,\sigma ^{-4}\right)
^{r/2}\,t_{1}^{2}\,\mbox{Var}(\delta _{i,n})\right.  \nonumber \\
&&\left. +\left( q_{i,n}+n^{-1}\right) ^{r/2-1}t_{2}^{2}\,\mbox{Var}%
(u_{i,n})\right)  \nonumber \\
&\leq &A\,r!\,\gamma _{n}^{r-2}\left( t_{1}^{2}\,\mbox{Var}(\delta
_{i,n})+t_{2}^{2}\,\mbox{Var}(u_{i,n})\right)  \nonumber \\
&\leq &A\,r!\,\gamma _{n}^{r-2}\,\mbox{Var}(t_{1}\,\delta
_{i,n}+t_{2}\,u_{i,n}),  \nonumber
\end{eqnarray}
for sufficiently large $n\geq n_{0}$. Using (\ref{gmm}) for $r=4$ and
H\"{o}lder's inequality, we get
\[
\left( \mbox{Var}(t_{1}\,\delta _{i,n}+t_{2}\,u_{i,n})\right) ^{2}\leq {\bf %
E\,}\left| t_{1\,}\delta _{i,n}+t_{2}\,u_{i,n}\right| ^{4}\leq A\,\gamma
_{n}^{2}\,\mbox{Var}(t_{1}\,\delta _{i,n}+t_{2}\,u_{i,n}).
\]
Hence,
\begin{equation}
\mbox{Var}(t_{1}\,\delta _{i,n}+t_{2}\,u_{i,n})\leq A\,\gamma _{n}^{2},\quad %
\mbox{for }\,\left\| t\right\| =1.  \label{hei0}
\end{equation}
Limit relation (\ref{fla}) follows from (\ref{tendzero}), (\ref{PSI}), (\ref
{gamma}), (\ref{qin}) and (\ref{777}).

We shall apply Lemma \ref{f3.1} with $m=s_{n}$,
\begin{equation}
H=A_{1},\quad \gamma =A_{2}\gamma _{n},\quad \lambda ^{2}=\min_{\left\|
t\right\| =1}\mbox{Var}(t_{1}\,S_{n}+t_{2}\,U_{n}),\quad \zeta _{i}=(\delta
_{i,n},u_{i,n}),  \label{hei1}
\end{equation}
\begin{equation}
\Lambda _{i}^{2}=\max_{\left\| t\right\| =1}\mbox{Var}(t_{1}\,\delta
_{i,n}+t_{2}\,u_{i,n})\leq 2\,\mbox{Var}(\delta _{i,n})+2\,\mbox{Var}%
(u_{i,n}),\quad i=1,\ldots ,s_{n},  \label{hei2}
\end{equation}
\begin{equation}
H^{*}=280\,H\,\lambda ^{-2}\sum_{i=1}^{s_{n}}\Lambda _{i}^{2},\quad \Xi
=(S_{n},U_{n})\in {\bf R}^{2},\quad {\it \ }\Theta ={\bf B}^{-1/2}\Xi ,
\label{hei3}
\end{equation}
where ${\bf B}$ is the covariance operator of $\Xi $. Fixing $A_{1}=A$ from (%
\ref{Bern}), using (\ref{hei0}) and (\ref{hei2}) and choosing $A_{2}$ to be
large enough, we ensure the validity of the inequality
\begin{equation}
18\,H\,\max_{1\leq i\leq s_{n}}\Lambda _{i}^{2}\leq \gamma ^{2}.
\label{hryu}
\end{equation}
Using (\ref{cn}), (\ref{cunion}), (\ref{p99}), (\ref{d0}), (\ref{vard}) and
Lemma \ref{l3.1}, we obtain (for sufficiently large $n\geq n_{0}$)
\begin{equation}
\sum_{i=1}^{s_{n}}\mbox{Var}(\delta _{i,n})=\sum_{i\in \Upsilon _{3}}\frac{%
\sigma _{n}^{2}(I_{i,n})}{\sigma _{n}^{2}(C_{n})}\leq 4\sum_{i\in \Upsilon
_{3}}\frac{{\bf P}(I_{i,n})\,\sigma ^{2}}{\sigma ^{2}}\leq 4.  \label{hei4}
\end{equation}
By (\ref{varuu}) and (\ref{vars}),
\begin{equation}
\sum_{i=1}^{s_{n}}\mbox{Var}(u_{i,n})=1-\alpha _{n}.  \label{hei5}
\end{equation}
Now (\ref{hei2}), (\ref{hei4}) and (\ref{hei5}) imply
\begin{equation}
\sum_{i=1}^{s_{n}}\Lambda _{i}^{2}\leq 10.  \label{hei6}
\end{equation}
Furthermore, by (\ref{tendzero}), (\ref{hn2}), (\ref{vars}), (\ref{hei1})
and inequality (\ref{cov}) of Lemma \ref{l3.4},
\begin{equation}
\lambda \geq \min \left\{ \mbox{Var}(S_{n}),\mbox{Var}(U_{n})\right\}
-2\,\left| \,\mbox{{\rm cov}}(S_{n},U_{n})\right| \geq 1/2,  \label{hei7}
\end{equation}
\begin{equation}
\mu \leq \max \left\{ \mbox{Var}(S_{n}),\mbox{Var}(U_{n})\right\} +2\,\left|
\,\mbox{{\rm cov}}(S_{n},U_{n})\right| \leq 2,  \label{hei10}
\end{equation}
for sufficiently large $n\geq n_{0}$, where $\mu $ is the maximal eigenvalue
of the covariance matrix ${\bf B}$. Applying Lemma \ref{f3.1} and taking
into account relations $\Xi ={\bf B}^{1/2}\Theta $, (\ref{aG}), (\ref{hei1}%
)--(\ref{hryu}) and (\ref{hei6})--(\ref{hei10}), we obtain, for $r\geq 3$, $%
n\geq n_{0}$:
\begin{eqnarray*}
\sup_{\left\| t\right\| =1}\left| \Gamma _{r}\left\{
t_{1}\,S_{n}+t_{2}\,U_{n}\right\} \right| &\leq &A^{r}\,\sup_{\left\|
t\right\| =1}\left| \Gamma _{r}\left\{ \left\langle t,\Theta \right\rangle
\right\} \right| \\
&\leq &A^{r}\,H^{*}\left( r-2\right) !\left( 8\sqrt{2}\,\gamma /\lambda
\right) ^{r-2}\leq \left( r-2\right) !\left( A\,\gamma _{n}\right) ^{r-2},
\end{eqnarray*}
proving (\ref{hei8}).\medskip

The following fact is well known. It may be easily derived from Remark 2 in
Rivlin \cite{r}, p. 96. It allows us to estimate coefficients of a
polynomial via its maximum on an interval.\medskip

\begin{lemma}
\label{f3.2} {\it Let }${\Bbb P}(x)=a_{0}+a_{1}x+\cdot \cdot \cdot
+a_{r}x^{r}$ {\it be a polynomial of degree not exceeding~}$r${\it . Then}
\[
\left| a_{k}\right| \leq \max \left\{ \left| t_{k}^{(r)}\right| ,\left|
t_{k}^{(r-1)}\right| \right\} \max_{-1\leq x\leq 1}\left| {\Bbb P}(x)\right|
,
\]
{\it where }$t_{k}^{(r)}$ {\it are coefficients of }${\Bbb
T}_{r}$, {\it the Chebyshev polynomial of order }$r$.\medskip
\end{lemma}

The Chebyshev polynomial
\[
{\Bbb T}_{r}(x)=t_{0}^{(r)}+t_{1}^{(r)}x+\cdot \cdot \cdot
+t_{r}^{(r)}x^{r},\quad r=1,2,\ldots ,
\]
is characterized as having the maximal leading coefficient $%
t_{r}^{(r)}=2^{r-1}$ among all polynomials ${\Bbb P}(x)$ with $\max_{-1\leq
x\leq 1}\left| {\Bbb P}(x)\right| \leq 1$. We have
\begin{equation}
{\Bbb T}_{0}(x)=1,\quad {\Bbb T}_{1}(x)=x,\quad {\Bbb T}_{r}(x)=2x\,{\Bbb T}%
_{r-1}(x)-{\Bbb T}_{r-2}(x),\quad r=2,3,\ldots ,  \label{Trx}
\end{equation}
see Rivlin (\cite{r}, formulas (1.11), (1.101)). By induction in $r$, it is
easy to derive from (\ref{Trx}) the rough bound
\begin{equation}
\sum_{k=0}^{r}\left| t_{k}^{(r)}\right| \leq 3^{r-1},\quad r=1,2,\ldots .
\label{3r}
\end{equation}

Let us consider the definition and some useful properties of classes of $d$%
-dimensional distributions \thinspace ${\cal A}_{d}(\tau )$, \thinspace $%
\tau \geq 0$, \thinspace introduced in Zaitsev \cite{z86}, see as well
Zaitsev \cite{z96}, \cite{z98} and \cite{z01}. The class \thinspace ${\cal A}%
_{d}(\tau )$ \thinspace (with a fixed \thinspace $\tau \geq 0$) consists of $%
d$-dimensional distributions~$F$ \thinspace for which the function
\[
\varphi (z)=\varphi (F,z)=\log \int_{{\bf R}^{d}}e^{\left\langle
z,x\right\rangle }F\{dx\}\qquad (\varphi (0)=0)
\]
is defined and analytic for \thinspace $\left\| z\right\| \tau <1$, $z\in
{\bf C}^{d}$, \thinspace and
\[
\left| d_{u}d_{v}^{2}\,\varphi (z)\right| \leq \Vert u\Vert \tau
\,\left\langle {\Bbb D}\,v,v\right\rangle \>\qquad \mbox{for all }\,u,v\in
{\bf R}^{d}\ \,\mbox{and}\ \,\left\| z\right\| \tau <1,
\]
where ${\Bbb D}$ is the covariance operator corresponding to $F$,\thinspace
and $d_{u}\varphi $ denotes the derivative\thinspace \thinspace of the
function $\varphi $ in direction $u$. It is easy to see that \thinspace $%
\tau _{1}<\tau _{2}$ implies \thinspace ${\cal A}_{d}\left( {\tau _{1}}%
\right) \subset {\cal A}_{d}\left( {\tau _{2}}\right) $. \thinspace
Moreover, the class ${\cal A}_{d}\left( {\tau }\right) $ is closed with
respect to convolution: if~\thinspace $F_{1},F_{2}\in {\cal A}_{d}\left( {%
\tau }\right) $, then~\thinspace $F_{1}*F_{2}\in {\cal A}_{d}\left( {\tau }%
\right) $. \thinspace The class ${\cal A}_{d}\left( {0}\right) $ coincides
with the class of all Gaussian distributions in ${\bf R}^{d}$.\medskip

\begin{lemma}
\label{l3.6} {\it For sufficiently large }$n\geq n_{0}${\it ,} {\it we have}
\begin{equation}
G\stackrel{\rm def}{=}{\cal L}\left( (S_{n},U_{n})\right) \in {\cal A}%
_{2}(\tau _{n}),\>\qquad \mbox{{\it where }}  \label{tau1}
\end{equation}
\begin{equation}
\tau _{n}=A\gamma _{n}=A\,\left( \Psi _{n}^{3/2}\,\,\max_{i\in \Upsilon
_{3}}\,p_{i,n}^{1/2}+\max_{1\leq i\leq s_{n}}\,q_{i,n}^{1/2}\right)
\rightarrow 0{\it ,\ as\ }n\rightarrow \infty ,  \label{tau}
\end{equation}
{\it with }$\Psi _{n}${\it \ defined in }{\rm (\ref{PSI})}.\medskip
\end{lemma}

\noindent {\it Proof.} Comparing formulas (\ref{cum1}) and (\ref{cum2}), we
see that
\begin{equation}
\frac{\Gamma _{r}\left\{ t_{1}\,S_{n}+t_{2}\,U_{n}\right\} }{r!}%
=\sum_{k=0}^{r}\frac{\Gamma _{k,r-k}\left\{ \left( S_{n},\,U_{n}\right)
\right\} \,t_{1}^{k}t_{2}^{r-k}}{k!\,(r-k)!},\quad r=1,2,\ldots .
\label{cumeq}
\end{equation}
Define polynomials ${\Bbb P}_{r}(x)=a_{0}^{(r)}+a_{1}^{(r)}x+\cdot \cdot
\cdot +a_{r}^{(r)}x^{r}$ with
\begin{equation}
a_{k}^{(r)}=\frac{\Gamma _{k,r-k}\left\{ \left( S_{n},\,U_{n}\right)
\right\} }{k!\,(r-k)!},\quad k=0,1,\ldots ,r.  \label{ak}
\end{equation}
By inequality (\ref{hei8}) of Lemma 3.5, (\ref{aG}) and (\ref{cumeq}), we
have, for $r=3,4,\ldots $:
\begin{equation}
\max_{-1\leq x\leq 1}\left| {\Bbb P}_{r}(x)\right| \leq \frac{1}{r!}%
\sup_{\left\| t\right\| \leq \sqrt{2}}\left| \Gamma _{r}\left\{
t_{1}\,S_{n}+t_{2}\,U_{n}\right\} \right| \leq \frac{2^{r/2}\left(
r-2\right) !}{r!}\left( A\,\gamma _{n}\right) ^{r-2},  \label{pol}
\end{equation}
if $n\geq n_{0}$ is sufficiently large. Applying Lemma \ref{f3.2} and
relations (\ref{3r}), (\ref{ak}) and (\ref{pol}), we get
\begin{eqnarray}
\left| \Gamma _{k,r-k}\left\{ \left( S_{n},\,U_{n}\right) \right\} \right|
&\leq &\frac{3^{r-1}2^{r/2}\left( r-2\right) !\,k!\,(r-k)!}{r!}\left(
A\,\gamma _{n}\right) ^{r-2},  \nonumber \\
&\leq &\left( r-2\right) !\,\left( A\,\gamma _{n}\right) ^{r-2},\quad
k=0,1,\ldots ,r,\quad r=3,4,\ldots .  \label{cumineq}
\end{eqnarray}
Further, expanding, for $u=(u_{1},u_{2})\in {\bf R}^{2}$, $%
v=(v_{1},v_{2})\in {\bf R}^{2}$, $w=(w_{1},w_{2})\in {\bf C}^{2}$,
\[
u=u_{1}e_{1}+u_{2}e_{2},\quad v=v_{1}e_{1}+v_{2}e_{2},\quad
w=w_{1}e_{1}+w_{2}e_{2},
\]
and rewriting $\Gamma _{r_{1},r_{2}}\left\{ \left( S_{n},\,U_{n}\right)
\right\} $ as
\[
\Gamma _{r_{1},r_{2}}\stackrel{\rm def}{=}\Gamma _{r_{1},r_{2}}\left\{
\left( S_{n},\,U_{n}\right) \right\} =\left.
d_{e_{1}}^{r_{1}}d_{e_{2}}^{r_{2}}\log {\bf E\,}\exp \left(
z_{1}\,S_{n}+z_{2}\,U_{n}\right) \right| _{z=0},
\]
we have
\begin{eqnarray*}
&&\left. d_{u}d_{v}^{2}d_{w}^{r}\,\log {\bf E\,}\exp \left(
z_{1}\,S_{n}+z_{2}\,U_{n}\right) \right| _{z=0} \\
&=&\sum_{k=0}^{r}\frac{r!}{k!\,(r-k)!}\,w_{1}^{k}w_{2}^{r-k}(\Gamma
_{k+3,r-k}\,u_{1}v_{1}^{2}+\Gamma
_{k+2,r+1-k}\,(u_{2}v_{1}^{2}+2\,u_{1}v_{1}v_{2}) \\
&&\qquad \qquad \qquad \qquad +\,\Gamma
_{k+1,r+2-k}\,(u_{1}v_{2}^{2}+2\,u_{2}v_{1}v_{2})+\Gamma
_{k,r+3-k}\,u_{2}v_{2}^{2}).
\end{eqnarray*}
Coupled with (\ref{cumineq}), this implies
\[
\left| \left. d_{u}d_{v}^{2}d_{w}^{r}\,\log {\bf E\,}\exp \left(
z_{1}\,S_{n}+z_{2}\,U_{n}\right) \right| _{z=0}\right| \leq r!\,\left\|
u\right\| \cdot \left\| v\right\| ^{2}\cdot \left\| w\right\| ^{r}\cdot
\left( A\,\gamma _{n}\right) ^{r+1},
\]
for $r=0,1,\ldots $. By Taylor's formula,
\[
\left. d_{u}d_{v}^{2}\,\log {\bf E\,}\exp \left(
z_{1}\,S_{n}+z_{2}\,U_{n}\right) \right| _{z=w}=\sum_{r=0}^{\infty }\frac{%
\left. d_{u}d_{v}^{2}d_{w}^{r}\,\log {\bf E\,}\exp \left(
z_{1}\,S_{n}+z_{2}\,U_{n}\right) \right| _{z=0}}{r!}.
\]
Therefore,
\[
\left| d_{u}d_{v}^{2}\,\log {\bf E\,}\exp \left(
z_{1}\,S_{n}+z_{2}\,U_{n}\right) \right| \leq A\,\gamma _{n}\,\left\|
u\right\| \cdot \left\| v\right\| ^{2},\quad \mbox{for }\left\| z\right\|
\cdot A\,\gamma _{n}\leq 1,
\]
for a suitably chosen absolute constant $A$. It remains to note that, by (%
\ref{tendzero}), (\ref{hn2}), (\ref{vars}) and (\ref{cov}),
\[
\mbox{{\rm Var}}(v_{1}\,S_{n}+v_{2}\,U_{n})=v_{1}^{2}\mbox{{\rm Var}}%
(S_{n})+v_{2}^{2}\,\mbox{{\rm Var}}(U_{n})+2\,v_{1}v_{2}\,\mbox{{\rm cov}}%
(S_{n},U_{n})\geq \left\| v\right\| ^{2}/2,
\]
for sufficiently large $n\geq n_{0}$. Limit relation (\ref{tau}) is a
consequence of (\ref{fla}). \medskip

\section{\bf Exponential bound for the integral over an exceptional set}

The proof of the following Lemma \ref{l4.2} is similar to the proof of
Gin\'{e}, Mason and Zaitsev \cite{gmz}, Proposition 3.1.$\medskip $

\begin{lemma}
\label{l4.2} {\it Let }$B$ {\it be a Borel subset of }${\bf R}$,{\it \ }
\begin{equation}
\xi _{n}=\int_{B}\left( \Delta _{n}(x)-{\bf E}\,\Delta _{n}(x)\right) \,dx
\label{a}
\end{equation}
{\it where}
\begin{equation}
\Delta _{n}(x)=\sqrt{n}\left| f_{n}(x)-{\bf E\,}f_{n}(x)\right| =\frac{1}{%
h_{n}\sqrt{n}}\left| \,\sum_{i=1}^{n}\left\{ K\left( \frac{x-X_{i}}{h_{n}}%
\right) -{\bf E}\,K\left( \frac{x-X}{h_{n}}\right) \right\} \right| {\bf .}
\label{doloto}
\end{equation}
{\it Then }
\begin{equation}
{\bf E}\,\exp \left\{ \lambda |\xi _{n}|\right\} \le 4\,\exp \left\{
\sum_{m=2}^{\infty }\left( \frac{720\,e\lambda \kappa }{\log m}\right)
^{m}\left( \Omega ^{m/2}(n,B)+\frac{1}{n^{m/2-1}}\,\Omega (n,B)\right)
\right\} ,  \label{14}
\end{equation}
{\it for all $\lambda \ge 0$.}$\medskip $ \end{lemma}

{\it Proof.} Let $X,X_{1},X_{1}^{\prime },X_{2},X_{2}^{\prime },\dots ,$ be
i.i.d. random variables. Further, we let $\eta $ be a Poisson random
variable with mean $n,$ independent of $X_{1},X_{1}^{\prime
},X_{2},X_{2}^{\prime },\dots ,$ and set
\[
\Delta _{\eta }(x)=\frac{1}{h_{n}\sqrt{n}}\left| \,\sum_{i=1}^{\eta }K\left(
\frac{x-X_{i}}{h_{n}}\right) -n\,{\bf E}\,K\left( \frac{x-X}{h_{n}}\right)
\right| .
\]
Define
\begin{equation}
\overline{\xi }_{n}=\int_{B}\left( \Delta _{n}(x)-{\bf E}\,\Delta _{\eta
}(x)\right) \,dx.  \label{b}
\end{equation}
Let ${\cal I}_{s}$, $s=1,\dots ,6$, be a partition of the integers ${\bf Z}$
such that:

i) if $i\ne j\in{\cal I}_s$ then $|i-j|\ge2$, and

ii) for every $s=1,\dots ,6$, $\sum_{i\in {\cal I}_{s}}{\bf P}\left\{ X\in
\left( (i-1/2)h_{n},(i+3/2)h_{n}\right] \right\} \le 1/2$,

\noindent and set
\[
A_{s}=\cup _{i\in {\cal I}_{s}}B_{j,n},\quad s=1,\dots ,6,\quad \mbox{where }%
B_{j,n}=\left( ih_{n},(i+1)h_{n}\right] \cap B,
\]
Now, replacing $K_{1}$, $K_{2}$, $\eta _{n}$ and $\left(
ih_{n},(i+1)h_{n}\right] $ in the proof of inequalities (3.5), (3.7), (3.8)
and (3.13) in Gin\'{e}, Mason and Zaitsev \cite{gmz} by $K$, $0$, $\eta $
and $B_{j,n}$, respectively, and using the arguments therein, we obtain
\begin{eqnarray}
\hspace{0.2cm}{\bf E}\,\exp \{\lambda |\xi _{n}|\} &\le &{\bf E}\,\exp
\left\{ 2\lambda \overline{|\xi }_{n}|\right\}  \label{16} \\
&\le &\prod_{s=1}^{6}\left( {\bf E}\,\exp \left\{ 12\lambda \left|
\int_{A_{s}}(\Delta _{n}(x)-{\bf E}\,\Delta _{\eta }(x))\,dx\,\right|
\right\} \right) ^{1/6}  \nonumber \\
&\le &2\,\prod_{s=1}^{6}\left( {\bf E}\,\exp \left\{ 12\lambda \left|
\int_{A_{s}}(\Delta _{\eta }(x)-{\bf E}\,\Delta _{\eta }(x))\,dx\,\right|
\right\} \right) ^{1/6}  \nonumber
\end{eqnarray}
and
\begin{eqnarray}
&&\exp \left\{ 12\lambda \left| \int_{A_{s}}\left( \Delta _{\eta }(x)-{\bf E}%
\,\Delta _{\eta }(x)\right) \,dx\,\right| \right\}  \label{19} \\
&\le &\,2\,\exp \Biggl\{\sum_{j\in {\cal I}_{s}}\sum_{m=2}^{\infty }\left(
\frac{720e\lambda }{\log m}\right) ^{m}\Biggl[ \left( \int_{B_{j,n}}\frac{1}{%
h_{n}}{\bf E}\,K^{2}\left( \frac{x-X}{h_{n}}\right) \,dx\right) ^{m/2}
\nonumber \\
&&+\,\frac{1}{n^{m/2-1}}\int_{B_{j,n}}\frac{1}{h_{n}}{\bf E}\,\left|
\,K\left( \frac{x-X}{h_{n}}\right) \right| ^{m}\,dx\Biggr]\Biggr\}.
\nonumber
\end{eqnarray}
Furthermore, by a change of variables,
\begin{eqnarray*}
\sum_{j\in {\cal I}_{s}}\!\left( \int_{B_{j,n}}\frac{1}{h_{n}}{\bf E}%
\,K^{2}\left( \frac{x-X}{h_{n}}\right) \,dx\right) ^{\!m/2}\hspace{-0.5cm}
&\le &\!\left( \sum_{j\in {\cal I}_{s}}\int_{B_{j,n}}\frac{1}{h_{n}}{\bf E}%
\,K^{2}\left( \frac{x-X}{h_{n}}\right) \,dx\right) ^{\!m/2} \\
&\le &\!\left( {\bf E}\,\int_{B}\frac{1}{h_{n}}K^{2}\left( \frac{x-X}{h_{n}}%
\right) \,dx\right) ^{m/2}
\end{eqnarray*}
Using (\ref{k1}), (\ref{k2}), (\ref{13}), (\ref{19b}) and (\ref{19x}), we
obtain
\begin{eqnarray*}
{\bf E}\,\frac{1}{h_{n}}\int_{B}K^{2}\left( \frac{x-X}{h_{n}}\right) \,dx
&\leq &\kappa ^{2}\,h_{n}^{-1}\int_{B}{\bf P}\{X\in
[x-h_{n}/2,x+h_{n}/2]\}\,dx \\
&\leq &\kappa ^{2}\,\Omega (n,B).
\end{eqnarray*}
Similarly, we have
\begin{eqnarray*}
&&\frac{1}{n^{m/2-1}}\sum_{j\in {\cal I}_{s}}\int_{B_{j,n}}\frac{1}{h_{n}}%
{\bf E}\,\left| \,K\left( \frac{x-X}{h_{n}}\right) \right| ^{m}\,dx \\
&\le &\frac{\kappa ^{m}\,h_{n}^{-1}}{n^{m/2-1}}\,\int_{B}{\bf P}\{X\in
[x-h_{n}/2,x+h_{n}/2]\}\,dx\le \frac{\kappa ^{m}}{n^{m/2-1}}\,\Omega (n,B).
\end{eqnarray*}
Then, combining these estimates with (\ref{13}) and (\ref{19}), we obtain
\begin{eqnarray}
&&{\bf E}\,\exp \left\{ 12\lambda \left| \int_{A_{s}}\left( \Delta _{\eta
}(x)-{\bf E}\,\Delta _{\eta }(x)\right) \,dx\right| \right\}  \label{ooo} \\
&\le &2\exp \left\{ \sum_{m=2}^{\infty }\left( \frac{720e\lambda \kappa }{%
\log m}\right) ^{m}\left( \Omega ^{m/2}(n,B)+\frac{1}{n^{m/2-1}}\,\Omega
(n,B)\right) \right\} .  \nonumber
\end{eqnarray}
Inequalities (\ref{16}) and (\ref{ooo}) imply (\ref{14}).\medskip

\section{{\bf Proof of Theorems \ref{t1.2}, \ref{t1.3} and \ref{t1.4}}}

Note now that for any absolute constant $A$ we have
\begin{equation}
A/\sqrt{n}\leq \tau _{n},  \label{antn}
\end{equation}
for sufficiently large $n\geq n_{0}$ (see (\ref{nq}) and (\ref{tau})).
Therefore, by Example 1.2 in Zaitsev \cite{z86},
\begin{equation}
H\stackrel{\rm def}{=}{\cal L}\left( \left( 0,V_{n}\right) \right) \in {\cal %
A}_{2}\left( A/\sqrt{n}\right) \subset {\cal A}_{2}\left( \tau _{n}\right) .
\label{A2t1}
\end{equation}
Hence, by (\ref{tau1}) and (\ref{A2t1}),
\begin{equation}
Q\stackrel{\rm def}{=}{\cal L}\left( \left( S_{n},U_{n}\right) +\left(
0,V_{n}\right) \right) \in {\cal A}_{2}\left( \tau _{n}\right)  \label{uhh}
\end{equation}
(recall that $(S_{n},U_{n})$ is independent of $V_{n}$).

The following below Lemmas \ref{f5.1} and \ref{f5.2} are proved in Zaitsev
\cite{z02}. They provide estimates of the rate of convergence in a lemma of
Beirlant and Mason \cite{bm}, see as well Gin\'{e}, Mason and Zaitsev \cite
{gmz}, Lemma 2.4.\medskip

\begin{lemma}
\label{f5.1} {\it Let }$(${\it for each }$n\in {\bf N)}${\it \ }$\eta _{1,n}$%
{\it \ and }$\eta _{2,n}${\it \ be independent Poisson random variables with
}$\eta _{1,n}${\it \ being Poisson }$(n(1-\alpha _{n}))${\it \ and }$\eta
_{2,n}${\it \ being Poisson }$(n\alpha _{n})${\it \ where }$\alpha _{n}\in
(0,1).${\it \ Denote }$\eta _{n}=\eta _{1,n}+\eta _{2,n}${\it \ and set }
\[
U_{n}=\frac{\eta _{1,n}-n(1-\alpha _{n})}{\sqrt{n}}\quad {\it \mbox{ and }}%
\quad V_{n}=\frac{\eta _{2,n}-n\alpha _{n}}{\sqrt{n}}.
\]
{\it Let }$\{S_{n}\}_{n=1}^{\infty }${\it \ be a sequence of random
variables such that for each }$n\in {\bf N}${\it , the random vector }$%
(S_{n},U_{n})${\it \ is independent of }$V_{n}${\it . Assume that }$%
\mbox{\rm Var}(S_{n})=1$,
\begin{equation}
{\cal L}\left( \left( S_{n},U_{n}\right) +\left( 0,V_{n}\right) \right) \in
{\cal A}_{2}\left( \tau _{n}\right) ,  \label{Q}
\end{equation}
{\it and}
\begin{equation}
\left| \chi _{n}\right| \leq 1/2,  \label{xicond}
\end{equation}
{\it where}
\begin{equation}
\chi _{n}=\mbox{\rm cov}\left( S_{n},U_{n}\right) .  \label{xi}
\end{equation}
{\it Then there exist absolute constants }$A_{3},A_{4},A_{5},A_{6}$ {\it %
such that, for }$\tau _{n}${\it \ satisfying the estimates}
\begin{equation}
5\alpha _{n}^{-1}\,\exp \left\{ -5\alpha _{n}/432\,\tau _{n}^{2}\right\}
\leq \tau _{n},  \label{ac}
\end{equation}
\begin{equation}
A_{3}\,n^{-1/2}\leq \tau _{n}\leq A_{4},  \label{crd}
\end{equation}
{\it and} {\it for any fixed }$n\in {\bf N}${\it \ and} $y>0${\it , one can
construct on a probability space random variables }$\zeta _{n}${\it \ and }$%
Z ${\it \ so that the distribution of }$\zeta _{n}$ {\it is the conditional
distribution of }$S_{n}${\it \ given }$\eta _{n}=n${\it , }$Z${\it \ is a
standard normal random variable and}
\begin{equation}
{\bf P}\left\{ \left| \sqrt{1-\chi _{n}^{2}}\,Z-\zeta _{n}\right| \geq
\,y\right\} \leq A_{5}\,\exp \left\{ -A_{6}\,y/\tau _{n}\right\} .
\label{pi2}
\end{equation}
 \end{lemma}\medskip

\begin{lemma}
\label{f5.2}{\bf \ }{\it Let the conditions of Lemma\/ }$\ref{f5.1}${\it \
be satisfied. Then there exists absolute constants }$%
A_{7},A_{8},A_{9},A_{10} $ {\it such that, for any fixed }$n\in {\bf N}$
{\it and\ }$b${\it \ satisfying }
\begin{equation}
A_{3}n^{-1/2}\leq \tau _{n}\leq A_{7}b,\qquad b\leq 1,  \label{brestr}
\end{equation}
{\it one can construct on a probability space random variables }$\zeta _{n}$%
{\it \ and }$Z${\it \ with distributions described in Lemma\/ }$\ref{f5.1}$
{\it so that, for any }$y>0${\it ,}
\begin{eqnarray}
&&{\bf P}\left\{ \left| \sqrt{1-\chi _{n}^{2}}\,Z-\zeta _{n}\right| \geq
A_{10}\,\exp \left\{ -b^{2}/72\tau _{n}^{2}\right\} +y\right\}  \label{t14}
\\
&\leq &A_{8}\,\exp \left\{ -A_{9}\,y/\tau _{n}\right\} +2\,{\bf P}\left\{
\left| \omega \right| >y/6\right\} ,  \nonumber
\end{eqnarray}
{\it where }$\omega ${\it \ have the centered normal distribution
with variance }$b^{2}${\it .}$\medskip $ \end{lemma}

Comparing Lemmas \ref{f5.1} and \ref{f5.2}, we observe that in Lemma \ref
{f5.1} the probability space depends essentially on $y$, while the statement
(\ref{t14}) of Lemma \ref{f5.2} is valid on the same probability space
(depending on $b$) for any $y>0$. However, (\ref{t14}) is weaker than (\ref
{pi2}) for some values of $y$. The same rate of approximation (as in (\ref
{pi2})) is contained in (\ref{t14}) if $b^{2}\geq 72\tau _{n}^{2}\log
(1/\tau _{n})$ and $y\geq b^{2}/\tau _{n}$ only.{\bf \medskip }

Now we return to the estimation and note that, for random variables $%
S_{n},U_{n}$ and $V_{n}$ defined in (\ref{A})--(\ref{Vn}), the conditions of
Lemmas \ref{f5.1} and \ref{f5.2} are satisfied (with $\tau _{n}$ defined in
{\bf (}\ref{tau}) and $\eta _{n}=\eta $) for $n\geq n_{0}$. Indeed, by (\ref
{pqin})--(\ref{qqp}), we have
\begin{equation}
\tau _{n}^{*}\geq \tau _{n},  \label{013}
\end{equation}
if the constants $A$ in {\bf (}\ref{taun*}) and {\bf (}\ref{tau}) are chosen
in a suitable way. Limit relation (\ref{taun*}) follows from (\ref{tendzero}%
), (\ref{PSI}) and (\ref{gamma}). By (\ref{taun*}), (\ref{aldef}) and (\ref
{013}), $\alpha _{n}$ is chosen so that condition (\ref{ac}) is satisfied
for $n\geq n_{0}$. Note that by (\ref{tau}) and (\ref{antn}), condition (\ref
{crd}) and the first inequality in (\ref{brestr}) are fulfilled for
sufficiently large $n\geq n_{0}$. Moreover, by (\ref{tendzero}), and (\ref
{cov}), $\chi _{n}$ (defined in (\ref{xi})) tends to zero, as $n\rightarrow
\infty $, and condition (\ref{xicond}) is satisfied for sufficiently large $%
n\geq n_{0}$. Thus, we can apply to $S_{n},U_{n},V_{n}$ the statements of
Lemmas \ref{f5.1} and \ref{f5.2}.

By Lemma \ref{f5.1}, for sufficiently large fixed $n\geq n_{0}$ and for any
fixed $y>0$, one can construct on a probability space random variables $%
\zeta _{n}$\ and $Z$\ so that the distribution of $\zeta _{n}$ is the
conditional distribution of $S_{n}$\ given $\eta =n$,
\begin{equation}
\zeta _{n}=_{d}\sigma _{n}^{-1}(C_{n})\int_{C_{n}}\left( \Delta _{n}(x)-{\bf %
E\,}\Delta _{\eta }(x)\right) \,dx  \label{002}
\end{equation}
(see (\ref{0q}), (\ref{deleta}), (\ref{A}) and (\ref{doloto})) and a
standard normal random variable $Z$\ so that
\begin{equation}
{\bf P}\left\{ \left| \sqrt{1-\chi _{n}^{2}}\,Z-\zeta _{n}\right| \geq
\,y\right\} \leq A_{5}\,\exp \left\{ -A_{6}\,y/\tau _{n}\right\} .
\label{004}
\end{equation}
By Lemma \ref{f5.2}, for sufficiently large fixed $n\geq n_{0}$ and for any
fixed $b$\ satisfying
\begin{equation}
\tau _{n}\leq A_{7}b,\qquad b\leq 1,  \label{001}
\end{equation}
one can construct on a probability space a random variable $\zeta _{n}$\
with distribution described in (\ref{002}) and a standard normal random
variable $Z$ so that, for any $y>0$,
\begin{eqnarray}
&&{\bf P}\left\{ \left| \sqrt{1-\chi _{n}^{2}}\,Z-\zeta _{n}\right| \geq
A_{10}\,\exp \left\{ -b^{2}/72\tau _{n}^{2}\right\} +y\right\}  \label{003}
\\
&\leq &A_{8}\,\exp \left\{ -A_{9}\,y/\tau _{n}\right\} +2\,{\bf P}\left\{
\left| \omega \right| >y/6\right\} ,  \nonumber
\end{eqnarray}
where $\omega $\ have the centered normal distribution with variance $b^{2}$.

In both cases described above we can apply Lemma A of Berkes and Philipp
\cite{bp} assuming that there exists a sequence of i.i.d. random variables $%
X_{1},$ $X_{2},\ldots $ with probability density $f$ and such that
\begin{equation}
\zeta _{n}=\sigma _{n}^{-1}(C_{n})\int_{C_{n}}\left( \Delta _{n}(x)-{\bf E\,}%
\Delta _{\eta }(x)\right) \,dx,  \label{DZETANX}
\end{equation}
where $\Delta _{n}(x)$ is defined in (\ref{doloto}).

By (\ref{fnx}), (\ref{fetax}), (\ref{tendzero}), (\ref{cn}), (\ref{deleta}),
(\ref{doloto}) and Lemma \ref{l2.4}, we have
\begin{eqnarray}
&&\int_{C_{n}}\left| {\bf E\,}\Delta _{n}(x)-{\bf E\,}\Delta _{\eta
}(x)\right| \,dx  \label{005} \\
&=&\sqrt{n}\int_{C_{n}}\left| {\bf E\,}|f_{n}(x)-{\bf \,E\,}f_{n}(x)|-{\bf %
E\,}|f_{\eta }(x)-{\bf \,E\,}f_{n}(x)|\right| \,dx  \nonumber \\
&\leq &\sqrt{n}\int_{E_{n}}\left| {\bf E\,}|f_{n}(x)-{\bf \,E\,}f_{n}(x)|-%
{\bf E\,}|f_{\eta }(x)-{\bf \,E\,}f_{n}(x)|\right| \,dx  \nonumber \\
&\leq &\frac{A\,\lambda (E_{n})\,\left\| K^{3}\right\| }{\left\|
K^{2}\right\| \sqrt{nh_{n}^{2}}}+\frac{A\,{\Bbb N}_{n}\sqrt{h_{n}}}{\sqrt{%
\left\| K^{2}\right\| }}\stackrel{\rm def}{=}y_{n}  \nonumber
\end{eqnarray}
and $y_{n}\rightarrow \infty $, as $n\rightarrow \infty $. Applying Lemma
\ref{l4.2} for $B=\overline{C}_{n}$, we see that
\begin{equation}
{\bf E}\,\exp \left\{ \lambda |\xi _{n}|\right\} \le 4\,\exp \left\{
\sum_{m=2}^{\infty }\left( \frac{720e\lambda \kappa }{\log m}\right)
^{m}\left( \Omega ^{m/2}(n,\overline{C}_{n})+\frac{1}{n^{m/2-1}}\,\Omega (n,%
\overline{C}_{n})\right) \right\} ,  \label{006}
\end{equation}
for all $\lambda \ge 0$, where
\begin{equation}
\xi _{n}=\int_{\overline{C}_{n}}\left( \Delta _{n}(x)-{\bf E}\,\Delta
_{n}(x)\right) \,dx.  \label{007}
\end{equation}
By (\ref{dn3}) and (\ref{Pn1}),
\begin{equation}
n^{1/2}\,\Omega _{n}\rightarrow \infty ,\qquad \mbox{as}\quad n\rightarrow
\infty ,  \label{008}
\end{equation}
since we assume $nh_{n}^{2}\rightarrow \infty $. Using (\ref{dn3}), (\ref
{dn2}), (\ref{006}) and (\ref{008}), we obtain that, for sufficiently large $%
n\geq n_{0}$,
\begin{equation}
{\bf E}\,\exp \left\{ \lambda |\xi _{n}|\right\} \le 4\,\exp \left\{
\sum_{m=2}^{\infty }\left( \frac{A_{11}\,\lambda \,\kappa \,\Omega _{n}^{1/2}%
}{\log m}\right) ^{m}\right\} ,  \label{009}
\end{equation}
for all $\lambda \ge 0$. It may be shown that there exists an absolute
constant $A$ such that
\[
\sum_{m=2}^{\infty }\left( \frac{\mu }{\log m}\right) ^{m}\leq A\,\exp
\left\{ \exp \left\{ A\,\mu \right\} \right\} ,\quad \quad \mbox{for all }%
\mu >0.
\]
Applying the exponential Chebyshev inequality coupled with (\ref{009}),
where $$\lambda =A_{12}\,\kappa ^{-1}\,\Omega _{n}^{-1/2}\,\log ^{*}\log
^{*}(z/A_{11}\,\kappa \,\Omega _{n}^{1/2})$$ and $A_{12}$ is sufficiently
small, we obtain that
\begin{equation}
{\bf P}\left\{ \left| \xi _{n}\right| \geq z\right\} \leq A\exp \left\{
-A^{-1}\,\kappa ^{-1}\,\Omega _{n}^{-1/2}z\,\log ^{*}\log ^{*}(z/A\,\kappa
\,\Omega _{n}^{1/2})\right\} ,\quad \mbox{for any }z>0.  \label{010}
\end{equation}

Inequalities (\ref{p99}), (\ref{004}), (\ref{005}) and (\ref{010}) imply
that, for any fixed $n\geq n_{0}$ and for any fixed $x>0$, one can construct
on a probability space a sequence of i.i.d. random variables $X_{1},$ $%
X_{2},\ldots ,$ and a standard normal random variable $Z$\ so that
\begin{eqnarray}
&&{\bf P}\left\{ \left| \int_{-\infty }^{\infty }\left( \Delta _{n}(x)-{\bf E%
}\,\Delta _{n}(x)\right) \,dx-\sigma \,Z\,\right| \geq y_{n}+z+x\right\}
\label{012} \\
&\leq &A\,\Big( \exp \left\{ -A^{-1}\,\sigma ^{-1}x/\tau _{n}\right\} +\exp
\left\{ -A^{-1}\kappa ^{-1}\,\Omega _{n}^{-1/2}z\right\}  \nonumber \\
&+&{\bf P}\left\{ \left| \left( \sigma -\sigma _{n}(C_{n})\sqrt{1-\chi
_{n}^{2}}\right) Z\,\right| \geq z/2\right\} \Big) ,\quad \quad
\mbox{for
any }z>0.  \nonumber
\end{eqnarray}
Similarly, using (\ref{003}) instead of (\ref{004}), we establish that, for
any fixed $n\geq n_{0}$ and for any fixed $b$\ satisfying (\ref{001}), one
can construct on a probability space a sequence of i.i.d. random variables $%
X_{1},$ $X_{2},\ldots $ and a standard normal random variable $Z$\ so that
\[
{\bf P}\left\{ \left| \int_{-\infty }^{\infty }\left( \Delta _{n}(x)-{\bf E}
\,\Delta _{n}(x)\right) \,dx-\sigma \,Z\,\right| \geq A\,\sigma \,\exp
\left\{ -b^{2}/72\tau _{n}^{2}\right\} +y_{n}+z+x\right\}
\]
\[
\leq A\,\Big( \exp \left\{ -A^{-1}\,\sigma ^{-1}x/\tau _{n}\right\} +\exp
\{-A^{-1}\,\kappa ^{-1}\,\Omega _{n}^{-1/2}z\,\log ^{*}\log ^{*}(z/A\,\kappa
\,\Omega _{n}^{1/2})\}
\]
\begin{equation}
+\;{\bf P}\left\{ \left| \left( \sigma -\sigma _{n}(C_{n})\sqrt{1-\chi
_{n}^{2}}\right) Z\right| \geq z/2\right\}  \label{altern}
\end{equation}
\[
+\;{\bf P}\left\{ b\,\left| Z\right| >A^{-1}\,\sigma ^{-1}x\right\} \Big) %
,\quad \quad \quad \mbox{for any }x,z>0.
\]
Now, by (\ref{sikap}), (\ref{JnBK}), (\ref{sigk}), (\ref{psin0}), (\ref{dn3}%
), (\ref{L4a}), (\ref{cn}), (\ref{ghhg}), (\ref{cov}) and (\ref{xi}), we
have
\begin{eqnarray}
&&\left| \,\sigma -\sigma _{n}(C_{n})\sqrt{1-\chi _{n}^{2}}\,\right|
\label{014} \\
&\leq &\left| \,\sigma -\sigma _{n}(C_{n})\right| +\sigma \,\left| \sqrt{%
1-\chi _{n}^{2}}-1\right|  \nonumber \\
&\leq &\sigma \,\left( 1-\sqrt{{\bf P}(E_{n})}\right) +\left| \,\sigma
_{n}(E_{n})-\sigma \,\sqrt{{\bf P}(E_{n})}\,\right| +\sigma
_{n}(E_{n}\backslash C_{n})+\sigma \,\chi _{n}^{2}  \nonumber \\
&\leq &\frac{A\,||K^{2}||}{\sigma \,h_{n}}\,\left( {\Bbb L}_{n}+\frac{%
\varepsilon _{n}\,{\Bbb M}_{n}}{\left\| K^{2}\right\| }\right) +\,A\,\kappa
\,\Omega _{n}^{1/2}+\frac{A}{\sigma }\left( \frac{\left\| K^{3}\right\|
\,\lambda (E_{n})}{\left\| K^{2}\right\| \sqrt{nh_{n}^{2}}}\right) ^{2},
\nonumber
\end{eqnarray}
for sufficiently large $n\geq n_{0}$. Now inequality (\ref{015}) follows
from (\ref{doloto}), (\ref{013}), (\ref{012}) and (\ref{014}). Relations (%
\ref{tendzero}) and (\ref{hjk}) imply the limit relation in (\ref{ddn}). The
proof of Theorem \ref{t1.3} repeats that of Theorem \ref{t1.2}. The only
difference is that we apply (\ref{altern}) instead of (\ref{012}). \medskip

{\it Proof of Theorem }\ref{t1.4}. Without loss of generality, we assume $%
x\geq 1$. By Theorem \ref{t1.2}, for any $z>0$,
\begin{eqnarray*}
&&1-F(x)\leq 1-\Phi \left( x-2z-y_{n}/\sigma \right) +A\,\left( \exp
\{-A^{-1}\,z/\tau _{n}^{*}\}\right. \\
&&+\;\left. \exp \{-A^{-1}\,\kappa ^{-1}\,\Omega _{n}^{-1/2}\sigma \,z\,\log
^{*}\log ^{*}(\sigma \,z/A\,\kappa \,\Omega _{n}^{1/2})\}+{\bf P}\left\{
\left| \partial _{n}Z\right| \geq \sigma \,z/2\right\} \right)
\end{eqnarray*}
and
\begin{eqnarray*}
&&1-F(x)\geq 1-\Phi \left( x+2z+y_{n}/\sigma \right) -A\,\left( \exp
\{-A^{-1}\,z/\tau _{n}^{*}\}\right. \\
&&+\;\left. \exp \{-A^{-1}\,\kappa ^{-1}\,\Omega _{n}^{-1/2}\sigma \,z\,\log
^{*}\log ^{*}(\sigma \,z/A\,\kappa \,\Omega _{n}^{1/2})\}+{\bf P}\left\{
\left| \partial _{n}Z\right| \geq \sigma \,z/2\right\} \right) .
\end{eqnarray*}
Choosing here $z=\max \left\{ \sqrt{\tau _{n}^{*}x},\,\,\,\Omega _{n}^{1/4}%
\sqrt{x}\left( \log ^{*}\log ^{*}(1/\,\Omega _{n})\right) ^{-1/2},\sqrt{%
\partial _{n}}\right\} $ and using elementary properties of normal
distribution function, we get the result. \medskip

\subsection*{Acknowledgment}

A big part of the results of this paper was obtained by the author while
visiting the Delaware University. He would like to thank Professor David M.
Mason for his hospitality and for the statement of the problem. Many steps
of proofs of the present paper repeat corresponding steps of a detailed
proof of Theorem \ref{t1.1} (due to Mason) from the joint paper Gin\'{e},
Mason and Zaitsev \cite{gmz}. The author thanks Paul Eggermont, Evarist
Gin\'{e} and David Mason for useful discussions. The author is grateful to
the anonymous referee for careful reading the manuscript and useful
suggestions.

\end{document}